\documentclass[12pt]{article}
\usepackage{amsmath,amsfonts,amssymb,amsthm}
\input amssym.def
\begin{document}
\def \Z{\Bbb Z}
\def \C{\Bbb C}
\def \R{\Bbb R}
\def \Q{\Bbb Q}
\def \N{\Bbb N}
\def \k{{\bf k}}

\def \A{{\mathcal{A}}}
\def \D{{\mathcal{D}}}
\def \E{{\mathcal{E}}}
\def \E{{\mathcal{E}}}
\def \H{\mathcal{H}}
\def \S{{\mathcal{S}}}
\def \wt{{\rm wt}}
\def \tr{{\rm tr}}
\def \span{{\rm span}}
\def \Res{{\rm Res}}
\def \Der{{\rm Der}}
\def \End{{\rm End}}
\def \Ind {{\rm Ind}}
\def \Irr {{\rm Irr}}
\def \Aut{{\rm Aut}}
\def \GL{{\rm GL}}
\def \Hom{{\rm Hom}}
\def \mod{{\rm mod}}
\def \ann{{\rm Ann}}
\def \ad{{\rm ad}}
\def \rank{{\rm rank}\;}
\def \<{\langle}
\def \>{\rangle}

\def \g{{\frak{g}}}
\def \sl{{\frak{sl}}}
\def \gl{{\frak{gl}}}

\def \be{\begin{equation}\label}
\def \ee{\end{equation}}
\def \bex{\begin{example}\label}
\def \eex{\end{example}}
\def \bl{\begin{lem}\label}
\def \el{\end{lem}}
\def \bt{\begin{thm}\label}
\def \et{\end{thm}}
\def \bp{\begin{prop}\label}
\def \ep{\end{prop}}
\def \br{\begin{rem}\label}
\def \er{\end{rem}}
\def \bc{\begin{coro}\label}
\def \ec{\end{coro}}
\def \bd{\begin{de}\label}
\def \ed{\end{de}}

\newcommand{\nno}{\nonumber}
\newcommand{\nord}{\mbox{\scriptsize ${\circ\atop\circ}$}}
\newtheorem{thm}{Theorem}[section]
\newtheorem{prop}[thm]{Proposition}
\newtheorem{coro}[thm]{Corollary}
\newtheorem{conj}[thm]{Conjecture}
\newtheorem{example}[thm]{Example}
\newtheorem{lem}[thm]{Lemma}
\newtheorem{rem}[thm]{Remark}
\newtheorem{de}[thm]{Definition}
\newtheorem{hy}[thm]{Hypothesis}
\makeatletter \@addtoreset{equation}{section}
\def\theequation{\thesection.\arabic{equation}}
\makeatother \makeatletter

\begin{center}
{\Large \bf Vertex algebras associated with elliptic affine Lie
algebras}
\end{center}

\begin{center}
{Haisheng Li and Jiancai Sun\\
Department of Mathematical Sciences\\
Rutgers University, Camden, NJ 08102}
\end{center}

\begin{abstract}
We associate elliptic affine Lie algebras with what are called
vertex $\C((z))$-algebras and their modules in a certain category.
In the course, we construct two families of Lie algebras closely
related to elliptic affine Lie algebras.
\end{abstract}

\section{Introduction}
Elliptic affine Lie algebras, similar to affine Lie algebras, are a
family of infinite-dimensional Lie algebras associated with
finite-dimensional simple Lie algebras. Both elliptic affine Lie
algebras and affine Lie algebras are special examples of general
Krichever-Novikov algebras (see \cite{kn1}, \cite{kn2}). Let $\g$ be
a finite-dimensional simple Lie algebra over $\C$. Associated to
$\g$, one has an (untwisted) affine Lie algebra $\hat{\g}=\g\otimes
\C[t,t^{-1}]\oplus \C {\bf k}$, which is the universal central
extension of the Lie algebra $\g\otimes \C[t,t^{-1}]$. Let $R$ be
the quotient algebra of the (commutative and associative) algebra
$\C[t,t^{-1},u]$ modulo the relation $u^{2}=t^{3}-2\beta t^{2}+t$,
where $\beta$ is a fixed complex number. The elliptic affine Lie
algebra associated to $\g$, which is alternatively denoted by
$\hat{\g}_{e}$ in this current paper, is the universal
(three-dimensional) central extension of the Lie algebra $\g\otimes
R$ (see \cite{brem}, \cite{bcf}). Elliptic affine Lie algebras are
different from ordinary affine Lie algebras in several ways. For
example, unlike that affine Lie algebras are naturally $\Z$-graded,
elliptic affine Lie algebras are only ``quasi-graded." Nevertheless,
a ``highest weight module" theory has been developed in \cite{she1},
\cite{she2}, and a free field realization for $\hat{\g}_{e}$ with
$\g=sl_{2}$ has been obtained in \cite{bcf}.

In this paper, we study elliptic affine Lie algebras in the context
of vertex algebras and their modules. As our main result, we
associate elliptic affine Lie algebras with what were called vertex
$\C((t))$-algebras in \cite{li-tqva} and their modules of a certain
type.

It has been long known (see \cite{b-va}, \cite{flm}; cf. \cite{fz},
\cite{dl}) that affine Lie algebras $\hat{\g}$ can be canonically
associated with vertex algebras. This association can be described
as follows: For any complex number $\ell$, let $\C_{\ell}$ denote
the $1$-dimensional $(\g[t]+ \C{\bf k})$-module with $\g[t]$ acting
trivially and with ${\bf k}$ acting as scalar $\ell$. Let
$V_{\hat{\g}}(\ell,0)$ denote the $\hat{\g}$-module induced from
$(\g[t]+\C{\bf k})$-module $\C_{\ell}$.  Set ${\bf 1}=1\otimes 1\in
V_{\hat{\g}}(\ell,0)$. We have that
$V_{\hat{\g}}(\ell,0)=U(\hat{\g}){\bf 1}$ and $\g[t]{\bf 1}=0$.
Furthermore, there exists an operator $d$ on $V_{\hat{\g}}(\ell,0)$
such that $d{\bf 1}=0$ and
$$[d,a(x)]=\frac{d}{dx}a(x)\ \ \ \mbox{ for }a\in \g,$$
where $a(x)=\sum_{n\in \Z}a(n)x^{-n-1}$. This $\hat{\g}$-module
$V_{\hat{\g}}(\ell,0)$ is often referred in literature as the
(universal) level-$\ell$ vacuum module.  The fact is (see \cite{fz},
cf. \cite{ll}) that there exists a unique vertex algebra structure
on $V_{\hat{\g}}(\ell,0)$ with ${\bf 1}$ as the vacuum vector and
with $Y(a(-1){\bf 1},x)=a(x)$ for $a\in \g$. Furthermore, for any
restricted $\hat{\g}$-module $W$ of level $\ell$, there exists a
unique $V_{\hat{\g}}(\ell,0)$-module structure $Y_{W}:
V_{\hat{\g}}(\ell,0)\rightarrow \Hom (W,W((x)))$, such that
$Y_{W}(a(-1){\bf 1},x)=a(x)$ for $a\in \g$. On the other hand, for
any $V_{\hat{\g}}(\ell,0)$-module $(W,Y_{W})$, $W$ is a restricted
$\hat{\g}$-module of level $\ell$ with
$$a(x)=Y_{W}(a(-1){\bf 1},x)\ \ \ \mbox{ for }a\in \g.$$
This correspondence provides an isomorphism between the category of
restricted $\hat{\g}$-modules of level $\ell$ and the category of
$V_{\hat{\g}}(\ell,0)$-modules.

As for elliptic affine Lie algebras, the situation is different in
an essential way. First of all, as it is shown in Section 3,
elliptic affine Lie algebras do {\em not} admit nontrivial (suitably
defined) vacuum modules. Then we have to find an alternative.
Motivated by the conceptual construction of vertex algebras and
their modules in \cite{li-local}, we consider restricted modules $W$
for elliptic affine Lie algebras in the sense that all the
generating functions of $\hat{\g}_{e}$ lie in $\Hom (W,W((x)))$.
(Note that elliptic affine Lie algebras do admit nontrivial
restricted modules.) For any restricted $\hat{\g}_{e}$-module $W$,
the generating functions form a local subset of $\Hom (W,W((x)))$,
just as with affine Lie algebra $\hat{\g}$, and it follows from
\cite{li-local} that the generating functions generate a vertex
algebra $V_{W}$ with $W$ as a canonical module. However, the vertex
algebra $V_{W}$ is not a $\hat{\g}_{e}$-module under the canonical
action, unlike the case for affine Lie algebra $\hat{\g}$. To a
certain extent, this phenomenon is similar to what we have
experienced in \cite{li-twisted} for twisted affine Lie algebras,
where the generating functions of a twisted affine Lie algebra
$\hat{\g}[\sigma]$ on a restricted module generate a vertex algebra
which under the canonical action is a module for the {\em untwisted}
affine Lie algebra $\hat{\g}$, but {\em not} for the twisted affine
Lie algebra $\hat{\g}[\sigma]$.

We then introduce another Lie algebra $\check{\g}_{e}$ over the
field $\C((z))$, where $z$ is a formal variable. To a certain
extent, Lie algebra $\check{\g}_{e}$ is a deformation of
$\hat{\g}_{e}$. The good thing about this new Lie algebra is that
$\check{\g}_{e}$ viewed as a Lie algebra over $\C$ admits vacuum
modules. For every complex number $\ell$, using induced module
construction we construct a universal vacuum module
$V_{\check{\g}_{e}}(\ell,0)$ and we prove that there exists a
canonical vertex algebra structure on this vacuum module. Though
$V_{\check{\g}_{e}}(\ell,0)$ is also a $\C((z))$-module, it is {\em
not} a vertex algebra over the field $\C((z))$. The vertex algebra
structure and the $\C((z))$-module structure on
$V_{\check{\g}_{e}}(\ell,0)$ are encoded into a structure of a
so-called vertex $\C((z))$-algebra.

The notion of vertex $\C((z))$-algebra is a special case of the
notion of weak quantum vertex $\C((z))$-algebra, which was
introduced in \cite{li-tqva}. A {\em vertex $\C((z))$-algebra} is
simply a vertex algebra $V$ over $\C$, equipped with a
$\C((z))$-module structure such that
$$Y(f(z)u,x)(g(z)v)=f(z+x)g(z)Y(u,x)v$$
for $f(z),g(z)\in \C((z)),\; u,v\in V$. (Note that as the map $Y$ is
not $\C((z))$-linear, a vertex $\C((z))$-algebra is {\em not} a
vertex algebra over the field $\C((z))$.) For a vertex
$\C((z))$-algebra $V$, we define (see \cite{li-tqva}) a {\em type
zero $V$-module} to be a module $(W,Y_{W})$ for $V$ viewed as a
vertex algebra over $\C$, satisfying
$$Y_{W}(f(z)v,x)=f(x)Y_{W}(v,x)\ \ \
\mbox{ for }f(z)\in \C((z)),\; v\in V.$$ As for elliptic affine Lie
algebras, we prove that a restricted $\hat{\g}_{e}$-module structure
of level $\ell$ on a vector space $W$ over $\C$ is equivalent to a
type zero $V_{\check{\g}_{e}}(\ell,0)$-module structure.

In fact, what we have done in this paper is more general. We start
with a (possibly infinite-dimensional) Lie algebra $\g$ over $\C$,
equipped with a non-degenerate symmetric invariant bilinear form
$\<\cdot,\cdot\>$. For any polynomial $p(x)\in \C[x]$, we construct
a Lie algebra $\hat{\g}_{p}$ with underlying vector space
$$\hat{\g}_{p}=(\g\oplus \g^{1})\otimes \C[t,t^{-1}]\oplus \C \k,$$
where $\g^{1}=\g$ as a vector space. When $p(x)=x^{3}-2\beta
x^{2}+x$, $\hat{\g}_{p}$ is isomorphic to the quotient algebra of
the elliptic affine Lie algebra $\hat{\g}_{e}$ modulo a
two-dimensional central ideal. We also construct another Lie algebra
$\check{\g}_{p}$ over $\C((z))$ with underlying vector space
$$\check{\g}_{p}=\C((z))\otimes (\g\oplus \g^{1})\otimes
\C[t,t^{-1}]\oplus \C((z))\k.$$ For any complex number $\ell$, we
construct a (universal) vacuum $\check{\g}_{p}$-module
$V_{\check{\g}_{p}}(\ell,0)$ of level $\ell$ and we show that there
exists a canonical vertex $\C((z))$-algebra structure on
$V_{\check{\g}_{p}}(\ell,0)$. Furthermore,  we establish a canonical
isomorphism between the category of type zero
$V_{\check{\g}_{p}}(\ell,0)$-modules and the category of restricted
$\hat{\g}_{p}$-modules of level $\ell$.

This paper is organized as follows: In Section 2, we recall from
\cite{li-tqva} the basic notions and results, including the
definition of a vertex $\C((z))$-algebra and that of a module for a
vertex $\C((z))$-algebra, and including the conceptual construction
of vertex $\C((z))$-algebras and their modules. We also study vertex
$\C((z))$-algebras associated to a certain Heisenberg Lie algebra.
In Section 3, we associate Lie algebras $\hat{\g}_{p}$ and
$\check{\g}_{p}$ to a Lie algebra $\g$ with a non-degenerate
symmetric invariant bilinear form and a polynomial $p(x)$, and we
construct a vertex $\C((z))$-algebra $V_{\check{\g}_{p}}(\ell,0)$
for every complex number $\ell$. We then establish an isomorphism
between the category of $\check{\g}_{p}$-modules of level $\ell$ and
the category of type zero $V_{\check{\g}_{p}}(\ell,0)$-modules.

\section{Vertex $\C((t))$-algebras and their modules}

In this section, we first recall from \cite{li-tqva} the notion of
vertex $\C((t))$-algebra and two categories of modules for a vertex
$\C((t))$-algebra, and we then study an infinite-dimensional
Heisenberg-type Lie algebra in the context of vertex
$\C((t))$-algebras and their modules.

In this paper, the scalar field will be the field $\C$ of complex
numbers, unless it is specified otherwise. We shall use the formal
variable notations and conventions as established in \cite{flm} (cf.
\cite{ll}). For any formal variable $t$, let $\C((t))$ denote the
ring of lower truncated formal Laurent series. In fact, $\C((t))$ is
a field. The following notion is a special case of the notion of
weak quantum vertex $\C((t))$-algebra which was introduced and
studied in \cite{li-tqva}:

\bd{dvtalgebra} {\em A {\em vertex $\C((t))$-algebra} is a vertex
algebra $V$ over $\C$, equipped with a $\C((t))$-module structure
such that
\begin{eqnarray}
Y(f(t)u,x)g(t)v=f(t+x)g(t)Y(u,x)v
\end{eqnarray}
for $f(t),g(t)\in \C((t)),\ u,v\in V$, where
$f(t+x)=e^{x\frac{d}{dt}}f(t)\in \C((t))[[x]]$.} \ed

Note that in Section 3, due to a notion conflict we shall use $z$
instead of $t$ and we then shall be dealing with vertex
$\C((z))$-algebras. For vertex $\C((t))$-algebras, the following two
categories of modules are of our interest.

\bd{dmodules} {\em Let $V$ be a vertex $\C((t))$-algebra. A {\em
type one $V$-module} is a module $(W,Y_{W})$ for $V$ viewed as a
vertex algebra over $\C$, equipped a $\C((t))$-module structure such
that
\begin{eqnarray}
Y_{W}(f(t)u,x)g(t)w=f(t+x)g(t)Y_{W}(u,x)w
\end{eqnarray}
for $f(t),g(t)\in \C((t)),\; u\in V,\; w\in W$. A {\em type zero
$V$-module} is a module $(W,Y_{W})$ for $V$ viewed as a vertex
algebra over $\C$ such that
\begin{eqnarray}
Y_{W}(f(t)u,x)w=f(x)Y_{W}(u,x)w
\end{eqnarray}
for $f(t)\in \C((t)),\; u\in V,\; w\in W$.} \ed

Let $W$ be a general vector space over $\C$. Set
$$\E(W)=\Hom (W,W((x)))\subset (\End W)[[x,x^{-1}]].$$
For $a(x),b(x)\in \E(W)$, we say $a(x)$ and $b(x)$ are {\em local}
if there exists a nonnegative integer $k$ such that
$$(x-z)^{k}a(x)b(z)=(x-z)^{k}b(z)a(x).$$ We say a subset
$U$ of $\E(W)$ is {\em local} if $a(x),b(x)$ are local for any
$a(x),b(x)\in U$. For $a(x),b(x)\in \E(W),\; n\in \Z$, define
$a(x)_{n}b(x)\in \E(W)$ by
\begin{eqnarray}\label{eanb-old}
a(x)_{n}b(x)
=\Res_{x_{1}}\left((x_{1}-x)^{n}a(x_{1})b(x)-(-x+x_{1})^{n}b(x)a(x_{1})\right).
\end{eqnarray}
It was proved in \cite{li-local} that any local subset $U$ of
$\E(W)$ generates a vertex algebra $\<U\>$, where the identity
operator $1_{W}$ on $W$ is the vacuum vector, and $W$ is a faithful
$\<U\>$-module with $Y_{W}(a(x),z)=a(z)$ for $a(x)\in \<U\>$.

\br{rdefinition} {\em For $a(x),b(x)\in \E(W)$, the definition of
$a(x)_{n}b(x)\in \E(W)$ for $n\in \Z$ was modified in
\cite{li-gamma} (cf. \cite{li-g1}). Since we shall use some results
of \cite{li-gamma}, we here recall the connection between the two
definitions. Assume
$$f(x,z)a(x)b(z)\in \Hom (W,W((x,z)))$$
for some nonzero polynomial $f(x,z)$. Then
$$\left(f(x_{1},x)a(x_{1})b(x)\right)|_{x_{1}=x+x_{0}}\in \Hom
(W,W((x))[[x_{0}]]).$$ Let $\iota_{x,x_{0}}$ denote the unique field
embedding of $\C(x,x_{0})$ into $\C((x))((x_{0}))$, extending the
identity endomorphism of $\C[x,x_{0}]$, where $\C(x,x_{0})$ denotes
the field of rational functions.  Then $a(x)_{n}b(x)\in \E(W)$ for
$n\in \Z$ were defined in \cite{li-gamma} in terms of generating
function
$$Y_{\E}(a(x),x_{0})b(x)=\sum_{n\in \Z}a(x)_{n}b(x)x_{0}^{-n-1}$$
by
\begin{eqnarray}
Y_{\E}(a(x),x_{0})b(x)
=\iota_{x,x_{0}}(1/f(x+x_{0},x))\left(f(x_{1},x)a(x_{1})b(x)\right)|_{x_{1}=x+x_{0}}.
\end{eqnarray}
It was shown therein that if $a(x),b(x)$ are local, then this
definition coincides with the definition (\ref{eanb-old}). However,
the two definitions give different objects in general.} \er

The following, which was proved in \cite{li-local}, is very useful
in determining the structure of vertex algebras generated by local
subsets:

\bp{pcommutator} Let $V$ be a vertex algebra, let $u, v, w^{(0)},
\dots, w^{(k)}\in V$, and let $(W, Y_{W})$ be a $V$-module. If
\begin{eqnarray*}
[Y(u,x_{1}),Y(v,x_{2})]
=\sum_{i=0}^{k}\frac{1}{i!}Y(w^{(i)},x_{2})\left(\frac{\partial}{\partial
x_{2}}\right)^{i}x_{2}^{-1}\delta\left(\frac{x_{1}}{x_{2}}\right)
\end{eqnarray*}
holds on $V$, which is equivalent to
\begin{eqnarray*}
u_{i}v=w^{(i)}\ \ \mbox{for }0\le i\le k\ \ \mbox{and }\ \ u_{i}v=0\
\ \mbox{for }i>k,
\end{eqnarray*}
then
\begin{eqnarray*}
[Y_{W}(u, x_{1}),Y_{W}(v,x_{2})]
=\sum_{i=0}^{k}\frac{1}{i!}Y_{W}(w^{(i)},x_{2})\left(\frac{\partial}{\partial
x_{2}}\right)^{i}x_{2}^{-1}\delta\left(\frac{x_{1}}{x_{2}}\right)
\end{eqnarray*}
holds on $W$. If $W$ is a faithful $V$-module, then the converse is
also true. \ep

Notice that for $f(x)\in \C((x)),\; a(x)\in \Hom (W,W((x)))$,
$$f(x)a(x)\in \Hom (W,W((x))).$$
Then $\Hom (W,W((x)))$ is naturally a $\C((x))$-module. We consider
$\E(W)$ as a $\C((t))$-module with
\begin{eqnarray}
f(t)a(x)=f(x)a(x)\ \ \ \ \mbox{ for }f(t)\in \C((t)),\; a(x)\in
\E(W).
\end{eqnarray}
A subspace $U$ of $\E(W)$ is said to be {\em closed} if
$$a(x)_{n}b(x)\in U\ \ \ \mbox{ for }a(x),b(x)\in U,\; n\in \Z.$$
We have (see \cite{li-tqva}):

\bp{pvta} Let $W$ be a vector space over $\C$. Any closed local
$\C((t))$-subspace $V$ of $\E(W)$, containing $1_{W}$, is a vertex
$\C((t))$-algebra, and $W$ is a faithful type zero $V$-module with
$Y_{W}(a(x),z)=a(z)$ for $a(x)\in V$. On the other hand, for any
local subset $U$ of $\E(W)$, $\C((t))\<U\>$ is a vertex
$\C((t))$-algebra with $W$ as a type zero module, where $\<U\>$ is
the vertex algebra (over $\C$) generated by $U$. \ep

The following follows immediately from \cite{li-local}:

\bp{pvta-module} Let $V$ be a vertex $\C((t))$-algebra, let $W$ be a
vector space over $\C$, and let $Y_{W}: V\rightarrow \Hom
(W,W((x)))$ be a $\C$-linear map. Set
$$U=\{ Y_{W}(v,x)\;|\; v\in V\}\subset \E(W).$$
Then $(W,Y_{W})$ carries the structure of a type zero $V$-module if
and only if $U$ is local and $Y_{W}$ is a homomorphism of vertex
$\C((t))$-algebras from $V$ to $\C((t))\<U\>$.\ep

Next, using the rest of this section we present a toy example of
vertex $\C((t))$-algebras.

\bd{dhf-algebra} {\em Let $f(z)\in \C[[z,z^{-1}]]$. We define an
infinite-dimensional Lie algebra ${\mathcal{H}}(f)$ over $\C$ with
generators $c, \ \beta_{n}$ for $n\in \Z$, subject to relations
\begin{eqnarray}\label{e2.5}
&&[c,{\mathcal{H}}(f)]=0,\nonumber\\
&&[\beta(z),\beta(w)]=\frac{1}{2}f'(w)z^{-1}\delta\left(\frac{w}{z}\right)c
+f(w)\frac{\partial}{\partial
w}z^{-1}\delta\left(\frac{w}{z}\right)c,
\end{eqnarray}
where $\beta(x)=\sum_{n\in \Z}\beta_{n}x^{-n-1}$ and
$z^{-1}\delta(w/z)=\sum_{n\in \Z}w^{n}z^{-n-1}$.} \ed

It is straightforward to show that $c$ and $\beta_{n}$ for $n\in \Z$
form a basis of ${\mathcal{H}}(f)$. We say that an
${\mathcal{H}}(f)$-module $W$ is of {\em level} $\ell\in \C$ if $c$
acts on $W$ as scalar $\ell$, and we say that an
${\mathcal{H}}(f)$-module $W$ is {\em restricted} if for any $w\in
W$, $\beta_{n}w=0$ for $n$ sufficiently large. Let $W$ be a
restricted ${\mathcal{H}}(f)$-module of level $\ell$. {}From the
defining relation (\ref{e2.5}), we have
$$(z-w)^{2}[\beta(z),\beta(w)]=0.$$
Then $\{ \beta(x)\}$ is a local subset of $\E(W)$. In view of
Proposition \ref{pvta}, $\beta(x)$ generates a vertex algebra
$V_{W}$. To determine the structure of $V_{W}$ completely we shall
need another Lie algebra.

\bd{dliealgebra-K} {\em Let $\ell\in \C$. We define a Lie algebra
$K(\ell)$ over $\C$ with a basis $\{\tilde{\beta}_{n},\;
\tilde{c}_{n}\;|\; n\in \Z\}$, where the Lie bracket relations are
given by
\begin{eqnarray}
&&[\tilde{c}_{n},K(\ell)]=0,\nonumber\\
&&[\tilde{\beta}_{m},\tilde{\beta}_{n}]=\frac{\ell}{2}(m-n)\tilde{c}_{m+n-1}
\ \ \ \mbox{ for }m,n\in \Z.
\end{eqnarray}
(It is straightforward to see that this indeed defines a Lie
algebra.)} \ed

Form generating functions
$$\tilde{\beta}(x)=\sum_{n\in \Z}\tilde{\beta}_{n}x^{-n-1},\ \
\ \ \tilde{c}(x)=\sum_{n\in \Z}\tilde{c}_{n}x^{-n-1}.$$ Then the
nontrivial bracket relations amount to
\begin{eqnarray}
[\tilde{\beta}(z),\tilde{\beta}(w)]
=\frac{\ell}{2}\tilde{c}'(w)z^{-1}\delta\left(\frac{w}{z}\right)
+\tilde{c}(w)\frac{\partial}{\partial
w}z^{-1}\delta\left(\frac{w}{z}\right)\ell,
\end{eqnarray}
where $\tilde{c}'(w)$ denotes the formal derivative of
$\tilde{c}(w)$.
 Set
$$K(\ell)_{+}
={\rm span}\{ \tilde{\beta}_{n},\; \tilde{c}_{n}\;|\; n\ge
0\}\subset K(\ell).$$ It is readily seen that $K(\ell)_{+}$ is a Lie
subalgebra. Viewing $\C$ as a trivial $K(\ell)_{+}$-module, we form
the induced $K(\ell)$-module
\begin{eqnarray}
V_{K(\ell)}=U(K(\ell))\otimes_{U(K(\ell)_{+})}\C.
\end{eqnarray}
Set ${\bf 1}=1\otimes 1\in V_{K(\ell)}$, and set
$$\tilde{\beta}=\tilde{\beta}_{-1}{\bf 1},\ \ \ \
\tilde{c}=\tilde{c}_{-1}{\bf 1}\in V_{K(\ell)}.$$ It is
straightforward to see that $K(\ell)$ admits a derivation $d$ such
that
\begin{eqnarray}
d(\tilde{\beta}_{n})=-n\tilde{\beta}_{n-1},\ \
d(\tilde{c}_{n})=-n\tilde{c}_{n-1}\ \ \ \mbox{ for }n\in \Z.
\end{eqnarray}
We see that $d$ preserves the subalgebra $K(\ell)_{+}$. Then it
follows that $d$ gives rise to a linear operator $\D$ on
$V_{K(\ell)}$ such that $\D{\bf 1}=0$ and
\begin{eqnarray}
[\D,\tilde{\beta}(x)]=\frac{d}{dx}\tilde{\beta}(x),\ \ \
[\D,\tilde{c}(x)]=\frac{d}{dx}\tilde{c}(x).
\end{eqnarray}
Therefore, by a theorem of Frenkel-Kac-Radul-Wang [FKRW] and
Meurman-Primc [MP], there exists a vertex algebra structure on
$V_{K(\ell)}$, which is uniquely determined by the condition that
${\bf 1}$ is the vacuum vector and
\begin{eqnarray}
Y(\tilde{\beta},x)=\tilde{\beta}(x),\ \ \ \
Y(\tilde{c},x)=\tilde{c}(x).
\end{eqnarray}
 As a vertex algebra,
$V_{K(\ell)}$ is generated by subset $\{\tilde{\beta}, \tilde{c}\}$,
and we have
\begin{eqnarray}
&&[Y(\tilde{c},z),Y(\tilde{c},w)]=0
=[Y(\tilde{c},z),Y(\tilde{\beta},w)],\nonumber\\
&&[Y(\tilde{\beta},z),Y(\tilde{\beta},w)]
=\frac{\ell}{2}Y(\D\tilde{c},w)z^{-1}\delta\left(\frac{w}{z}\right)
+Y(\tilde{c},w)\frac{\partial}{\partial
w}z^{-1}\delta\left(\frac{w}{z}\right)\ell.\ \ \ \ \ \ \ \
\end{eqnarray}

\bp{phf-modules} Let $\ell\in \C,\; f(z)\in \C((z))$, and let $W$ be
any restricted ${\mathcal{H}}(f)$-module of level $\ell$. Then there
exists a $V_{K(\ell)}$-module structure $Y_{W}$ on $W$, which is
uniquely determined by
$$Y_{W}(\tilde{\beta},x)=\beta(x),\ \ \  Y_{W}(\tilde{c},x)=f(x).$$
On the other hand, let $(W,Y_{W})$ be an irreducible
$V_{K(\ell)}$-module. Then $Y_{W}(\tilde{c},x)=f(x)$ for some
$f(x)\in \C((x))$ and $W$ is a restricted ${\mathcal{H}}(f)$-module
of level $\ell$ with $\beta(x)$ acting as $Y_{W}(\tilde{\beta},x)$.
\ep

\begin{proof} For the first assertion, the uniqueness is clear as
$V_{K(\ell)}$ is generated by $\{\tilde{\beta}, \tilde{c}\}$. As for
the existence, set
$$U=\{\beta(x), f(x)\}\subset \E(W).$$
{}From the defining relations of  ${\mathcal{H}}(f)$, $U$ is local.
By Proposition \ref{pvta}, we have a vertex algebra $\<U\>$ with $W$
as a faithful module. With the relations (\ref{e2.5}), by
Proposition \ref{pcommutator}, we have
\begin{eqnarray*}
&&[Y_{\E}(\beta(x),z),Y_{\E}(\beta(x),w)]\\
&=&\frac{\ell}{2}
Y_{\E}(f'(x),w)z^{-1}\delta\left(\frac{w}{z}\right)
+Y_{\E}(f(x),w)\frac{\partial}{\partial
w}z^{-1}\delta\left(\frac{w}{z}\right)\ell.
\end{eqnarray*}
We also have
\begin{eqnarray*}
Y_{\E}(f'(x),w)=f'(x+w)=\frac{\partial}{\partial
w}f(x+w)=\frac{\partial}{\partial w}Y_{\E}(f(x),w).
\end{eqnarray*}
It follows that $\<U\>$ is a $K(\ell)$-module with
$\tilde{\beta}(z)$ acting as $Y_{\E}(\beta(x),z)$, $\tilde{c}(z)$ as
$Y_{\E}(f(x),z)$. By the construction of $V_{K(\ell)}$, there exists
a $K(\ell)$-module homomorphism $\psi$ from $V_{K(\ell)}$ to $\<U\>$
with $\psi({\bf 1})=1_{W}$. We have
$$\psi(\tilde{\beta})=\psi(\tilde{\beta}_{-1}{\bf 1})=\beta(x)_{-1}1_{W}=\beta(x),\ \
\psi(\tilde{c})=f(x)_{-1}1_{W}=f(x).$$ As
$\{\tilde{\beta},\tilde{c}\}$ generates $V_{K(\ell)}$, it follows
that $\psi$ is a homomorphism of vertex algebras. Consequently, $W$
becomes a $V_{K(\ell)}$-module where
$$Y_{W}(\tilde{\beta},z)=Y_{W}(\beta(x),z)=\beta(z),\ \ \
Y_{W}(\tilde{c},z)=Y_{W}(f(x),z)=f(z).$$

On the other hand, let $W$ be an irreducible $V_{K(\ell)}$-module.
We have
\begin{eqnarray*}
&&[Y_{W}(\tilde{c},z),Y_{W}(\tilde{c},w)]=0=[Y_{W}(\tilde{c},z),Y_{W}(\tilde{\beta},w)]
\end{eqnarray*}
and
\begin{eqnarray*}
&&[Y_{W}(\tilde{\beta},z),Y_{W}(\tilde{\beta},w)]\\
&=&\frac{\ell}{2}Y_{W}(\D\tilde{c},w)z^{-1}\delta\left(\frac{w}{z}\right)
+Y_{W}(\tilde{c},w)\frac{\partial}{\partial
w}z^{-1}\delta\left(\frac{w}{z}\right)\ell\\
&=&\frac{\ell}{2}\left(\frac{\partial}{\partial
w}Y_{W}(\tilde{c},w)\right)z^{-1}\delta\left(\frac{w}{z}\right)
+Y_{W}(\tilde{c},w)\frac{\partial}{\partial
w}z^{-1}\delta\left(\frac{w}{z}\right)\ell.
\end{eqnarray*}
Then $W$ is a $K(\ell)$-module with $\tilde{\beta}(x)$ and
$\tilde{c}(x)$ acting as $Y_{W}(\tilde{\beta},x)$ and
$Y_{W}(\tilde{c},x)$, respectively. Since $W$ is an irreducible
$V_{K(\ell)}$-module, it follows that $W$ is an irreducible
$K(\ell)$-module. As $K(\ell)$ is of countable dimension over $\C$,
so is $W$. Because $\tilde{c}_{n}$ for $n\in \Z$ are central, in
view of Schur's Lemma, we have $\tilde{c}(x)=f(x)\in
\C[[x,x^{-1}]]$. As $\tilde{c}(x)=Y_{W}(\tilde{c},x)\in \Hom
(W,W((x)))$ {}from the definition of a module, we have $f(x)\in
\C((x))$. Therefore, $W$ is an ${\mathcal{H}}[f]$-module of level
$\ell$ with $\beta(x)$ acting as $Y_{W}(\tilde{\beta},x)$.
\end{proof}

Next, we give a refinement of Proposition \ref{phf-modules} in terms
of vertex $\C((t))$-algebras and their type zero modules. Equip
$\C((t))$ with the vertex algebra structure given by the Borcherds
construction with derivation $\frac{d}{dt}$, where
$$Y(p(t),x)q(t)=p(t+x)q(t)\ \ \ \mbox{ for }p(t),q(t)\in \C((t)).$$
Then the tensor product vertex algebra $\C((t))\otimes V_{K(\ell)}$
over $\C$ is naturally a vertex $\C((t))$-algebra with $\C((t))$
acting on the first factor.

\bd{dvf} {\em Let $f(t)\in \C((t))$. Define $V[f]$ to be the
quotient vertex algebra of $\C((t))\otimes V_{K(\ell)}$ over $\C$
modulo the relation
\begin{eqnarray}
f(t)\otimes {\bf 1}=1\otimes \tilde{c}.
\end{eqnarray}} \ed

{}From the construction of the vertex algebra $\C((t))$, for
$p(t)\in \C((t))$, $p(t)_{-1}$ (a component of the vertex operator
$Y(p(t),x)$) is the left multiplication by $p(t)$. It follows that
$V[f]$ is also a quotient $\C((t))$-module. Then $V[f]$ is naturally
a vertex $\C((t))$-algebra.

\bt{theisenberg} Let $\ell\in \C,\; f(z)\in \C((z))$. For any
level-$\ell$ restricted ${\mathcal{H}}(f)$-module $W$, there exists
a type zero $V[f]$-module structure $Y_{W}$ (on $W$) which is
uniquely determined by
$$Y_{W}(\tilde{\beta},x)=\beta(x),\ \ Y_{W}(\tilde{c},x)=f(x).$$
 On the other hand, for
any type zero $V[f]$-module $(W,Y_{W})$, $W$ becomes a level-$\ell$
restricted ${\mathcal{H}}(f)$-module with $\beta(x)$ acting as
$Y_{W}(\tilde{\beta},x)$. \et

\begin{proof} Let $W$ be a level-$\ell$ restricted
${\mathcal{H}}(f)$-module. By Proposition \ref{phf-modules}, there
exists a $V_{K(\ell)}$-module structure $Y_{W}$ on $W$ with the
desired properties. Extend $Y_{W}$ to a $\C$-linear map
$$\bar{Y}_{W}: \C((t))\otimes V_{K(\ell)}\rightarrow \Hom (W,W((x)))$$
by
$$\bar{Y}_{W}(p(t)\otimes v,x)=p(x)Y_{W}(v,x)$$
for $p(t)\in \C((t)), \; v\in V_{K(\ell)}$. Note that
$\C((x))\;(=\C((x))1_{W})\subset \E(W)$ and that
$$Y_{\E}(p(x),z)q(x)=p(x+z)q(x)\ \ \ \mbox{ for }p(x),q(x)\in
\C((x)).$$ Thus $\C((x))$ is a vertex algebra with $W$ as a module.
That is, $W$ is a $\C((t))$-module with $Y_{W}(p(t),x)=p(x)$ for
$p(t)\in \C((t))$. It follows that $(W,\bar{Y}_{W})$ carries the
structure of a module for the tensor product vertex algebra
$\C((t))\otimes V_{K(\ell)}$ (over $\C$). Furthermore, it is readily
seen that $(W,\bar{Y}_{W})$ is a type zero module. Since
$$\bar{Y}_{W}(f(t)\otimes {\bf
1},x)=f(x)=Y_{W}(\tilde{c},x)=\bar{Y}_{W}(1\otimes \tilde{c},x),$$
it follows that $\bar{Y}_{W}$ reduces to a module structure for
$V[f]$ viewed as a vertex algebra over $\C$. This makes $W$ a type
zero $V[f]$-module.

On the other hand, let $(W,Y_{W})$ be a type zero $V[f]$-module. As
$V_{K(\ell)}$ is a vertex subalgebra of $\C((t))\otimes V_{K(\ell)}$
and $V[f]$ is a quotient vertex algebra of $\C((t))\otimes
V_{K(\ell)}$, $W$ is naturally a $V_{K(\ell)}$-module. Furthermore,
we have
$$Y_{W}(\tilde{c},x)=Y_{W}(f(t),x)=f(x)$$
and {}from [DLM] (cf. [LL]) we have
$$Y_{W}(\D
\tilde{c},x)=\frac{d}{dx}Y_{W}(\tilde{c},x)=f'(x).$$ Using all these
facts and Proposition \ref{pcommutator} we obtain
\begin{eqnarray*}
&&[Y_{W}(\tilde{\beta},z),Y_{W}(\tilde{\beta},w)]\\
&=&\frac{\ell}{2}Y_{W}(\D
\tilde{c},w)z^{-1}\delta\left(\frac{w}{z}\right) +
Y_{W}(\tilde{c},w)\frac{\partial}{\partial
w}z^{-1}\delta\left(\frac{w}{z}\right)\ell\\
&=&\frac{\ell}{2}f'(w)z^{-1}\delta\left(\frac{w}{z}\right) +
f(w)\frac{\partial}{\partial
w}z^{-1}\delta\left(\frac{w}{z}\right)\ell.
\end{eqnarray*}
This implies that $W$ is a level-$\ell$ ${\mathcal{H}}(f)$-module
with $\beta(x)$ acting as $Y_{W}(\tilde{\beta},x)$.
\end{proof}

\section{Vertex $\C((z))$-algebras associated with elliptic affine Lie algebras}
In this section, we shall associate elliptic affine Lie algebras and
their restricted modules with vertex $\C((z))$-algebras and their
type zero modules. More generally, for each polynomial $p(x)\in
\C[x]$, we construct a Lie algebra $\hat{\g}_{p}$ over $\C$, which
generalizes elliptic affine Lie algebra $\hat{\g}_{e}$ in a certain
way, and we also construct a Lie algebra $\check{\g}_{p}$ over
$\C((z))$.  Then we construct a vertex $\C((z))$-algebra
$V_{\check{\g}_{p}}(\ell,0)$ associated with $\check{\g}_{p}$ and a
complex number $\ell$, and we establish an isomorphism between the
category of restricted $\hat{\g}_{p}$-modules of level $\ell$ and
the category of type zero $V_{\check{\g}_{p}}(\ell,0)$-modules.

We begin by recalling elliptic affine Lie algebras, following
\cite{bcf}. Let $\beta$ be a complex number which is fixed
throughout this section. Denote by $A_{\beta}[t^{\pm 1},u]$ the
quotient algebra of $\C[t,t^{-1},u]$ modulo relation
$$u^{2}=t^{3}-2\beta t^{2}+t.$$
Let $\g$ be a (possibly infinite-dimensional) Lie algebra over $\C$,
equipped with a non-degenerate symmetric invariant bilinear form
$\<\cdot,\cdot\>$. Roughly speaking, the elliptic affine Lie algebra
associated to $\g$, denoted by $\hat{\g}_{e}$, is the universal
central extension of the Lie algebra $\g\otimes A_{\beta}[t^{\pm
1},u]$. Following \cite{bcf}, let $\g^{1}$ be a vector space
isomorphic to $\g$, with a fixed linear isomorphism $a\in \g\mapsto
a^{1}\in \g^{1}$. Then
\begin{eqnarray}
\hat{\g}_{e}=(\g\oplus \g^{1})\otimes \C[t,t^{-1}]\oplus \C {\bf
k}\oplus \C {\bf k}_{+}\oplus \C {\bf k}_{-}
\end{eqnarray}
as a vector space over $\C$, where ${\bf k}, {\bf k}_{\pm}$ are
central and the nontrivial bracket relations are written in terms of
generating functions for $a\in \g$:
$$a(x)=\sum_{n\in \Z}(a\otimes t^{n})x^{-n-1}, \ \ \ \
a^{1}(x)=\sum_{n\in \Z}(a^{1}\otimes t^{n})x^{-n-1}.$$ The
nontrivial bracket relations are
\begin{eqnarray*}
&&[a(z),b(w)]=[a,b](w)z^{-1}\delta\left(\frac{w}{z}\right)
+\<a,b\>\frac{\partial}{\partial
w}z^{-1}\delta\left(\frac{w}{z}\right){\bf k},\\
&&[a(z),b^{1}(w)]=[a,b]^{1}(w)z^{-1}\delta\left(\frac{w}{z}\right)
+\<a,b\> A(w)\frac{\partial}{\partial w}z^{-1}\delta\left(\frac{w}{z}\right),\\
&&[a^{1}(z),b^{1}(w)]\\
&&\ \ \ \ =\left((w^{3}-2\beta
w^{2}+w)[a,b](w)+\frac{1}{2}\<a,b\>(3w^{2}-4\beta w+1){\bf k}\right)
z^{-1}\delta\left(\frac{w}{z}\right)\\
&&\hspace{1cm}  +\<a,b\>(w^{3}-2\beta
w^{2}+w)\frac{\partial}{\partial
w}z^{-1}\delta\left(\frac{w}{z}\right){\bf k},
\end{eqnarray*}
where
$$A(w)=w(P(w^{-1},\beta)+P(w,\beta)){\bf
k}_{+}+w(Q(w^{-1},\beta)+Q(w,\beta)-2) {\bf k}_{-},$$ in which
$P(x,\beta)$ and $Q(x,\beta)$ are certain nonnegative power series
in $x$, depending on $\beta$.

Note that we have slightly modified the original relation (by
replacing the invariant bilinear form $\<\cdot,\cdot\>$ on $\g$ to
$-\<\cdot,\cdot\>$), so that the elliptic affine Lie algebra
$\hat{\g}_{e}$ contains the standard affine Lie algebra $\hat{\g}$
as a subalgebra.

We shall be interested in $\hat{\g}_{e}$-modules on which central
elements ${\bf k}$ and ${\bf k}_{\pm}$ act as complex scalars $\ell$
and $\ell_{\pm}$.

\bd{dlevel}{\em If $W$ is a $\hat{\g}_{e}$-module on which ${\bf
k}$, ${\bf k}_{+}$ and ${\bf k}_{-}$ act as complex scalars $\ell,\;
\ell_{\pm}$, respectively, we say that $W$ is of \emph{level
$(\ell,\ell_{+},\ell_{-})$}.}\ed

Just as with the ordinary affine Lie algebra $\hat{\g}$, for any
$\hat{\g}_{e}$-module $W$ and for $a\in \g,\; n\in\Z$,  we write
$a(n), a^{1}(n)$ for the operators on $W$, corresponding to
$a\otimes t^{n}, a^{1}\otimes t^{n}$, respectively, and we view
$a(x)$ and $a^{1}(x)$ as elements of $(\End W)[[x,x^{-1}]]$.

\bd{drestricted}{\em A $\hat{\g}_{e}$-module $W$ is called a
\emph{restricted module} if for every $a\in \g$ and $w\in W$,
$a(n)w=a^{1}(n)w=0$ for $n$ sufficiently large, namely, if
$$a(x),a^{1}(x)\in \Hom (W,W((x)))\ \ \ \mbox{ for }a\in \g.$$
}\ed

We have (cf. \cite{bcf}, Lemma 4.3):

\bl{lkpm=0} Let $W$ be a restricted $\hat{\g}_{e}$-module on which
${\bf k}_{\pm}$ act as scalars $\ell_{\pm}$. Then
$$A(z)=z(P(z^{-1},\beta)+P(z,\beta))\ell_{+}
+z(Q(z^{-1},\beta)+Q(z,\beta)-2)\ell_{-}\in \C((z)).$$ \el

\begin{proof} Let $a,b\in \g$ with $\<a,b\>\ne 0$.
{}From the defining relations we get
\begin{eqnarray*}
\Res_{z}(z-w)[a(z),b^{1}(w)] =\<a,b\>A(w).
\end{eqnarray*}
Since $a(w),b^{1}(w)\in \Hom (W,W((w)))$, we have
$$\Res_{z}(z-w)[a(z),b^{1}(w)]\in \Hom (W,W((w))).$$
Thus $A(w)\in \C((w))$.
\end{proof}

{}From the defining relations we have
\begin{eqnarray*}
(z-w)^{2}[a(z),b(w)]=(z-w)^{2}[a^{1}(z),b^{1}(w)]
=(z-w)^{2}[a(z),b^{1}(w)]=0
\end{eqnarray*}
for $a,b\in \g$. Let $W$ be a restricted $\hat{\g}_{e}$-module. Set
$$U_{W}=\{a(x), a^{1}(x)\;|\; a\in \g\}\subset \E(W).$$
Then $U_{W}$ is a local subset. Thus $U_{W}$ generates a vertex
algebra $V_{W}$ inside $\E(W)$ with $W$ as a faithful module.

Note that if $W$ is a level-$\ell$ restricted module for the affine
Lie algebra $\hat{\g}$, the vertex algebra generated by the
generating functions $a(x)$ for $a\in \g$ is a vacuum
$\hat{\g}$-module of level $\ell$ with $a(z)$ acting as
$Y_{\E}(a(x),z)$. For elliptic affine Lie algebra $\hat{\g}_{e}$,
this is {\em no longer} the case. In fact, the following proposition
asserts that there does not exist a nontrivial vacuum
$\hat{\g}_{e}$-module.

\bp{pvacuum-module} Suppose that $W$ is a restricted
$\hat{\g}_{e}$-module of level $(\ell,\ell_{+},\ell_{-})$, equipped
with a vector $w_{0}\in W$ and a linear operator $D$ on $W$ such
that $W=U(\hat{\g}_{e})w_{0}$,
\begin{eqnarray*}
&&Dw_{0}=0,\ \ \ a(x)w_{0},\ a^{1}(x)w_{0}\in W[[x]],\\
 &&[D,a(x)]=\frac{d}{dx}a(x),\
\ \ [D,a^{1}(x)]=\frac{d}{dx}a^{1}(x)\ \ \ \mbox{ for }a\in \g.
\end{eqnarray*}
Then $\ell=0$ and $A(z)\in \C$. Furthermore, if $[\g,\g]=\g$, then
$W$ is a $1$-dimensional trivial $\hat{\g}_{e}$-module. \ep

\begin{proof} With $W$ assumed to be a restricted
$\hat{\g}_{e}$-module, $a(x),a^{1}(x)$ for $a\in \g$ form a local
subset $U_{W}$ of $\E(W)$. Then we have a vertex algebra $V_{W}$
generated by $U_{W}$. By Proposition 5.4.1 of \cite{ll}, we have
$$[D,\psi(x)]=\frac{d}{dx}\psi(x)\ \ \ \mbox{ for }\psi(x)\in
V_{W}.$$ Let $a,b\in \g$ with $\<a,b\>\ne 0$. We have
\begin{eqnarray*}
a(x)_{1}b^{1}(x)=\Res_{z}(z-x)[a(z),b^{1}(x)] =\<a,b\>A(x),
\end{eqnarray*}
which implies $A(x)\in V_{W}$. Then we have
$\frac{d}{dx}A(x)=[D,A(x)]=0$, proving $A(x)\in \C$. Similarly, we
have
\begin{eqnarray*}
a^{1}(x)_{1}b^{1}(x)=\Res_{x_{1}}(x_{1}-x)[a^{1}(x_{1}),b^{1}(x)]=\ell\<a,b\>(x^{3}-2\beta
x^{2}+x).
\end{eqnarray*}
Then
$$\ell\<a,b\>(3x^{2}-4\beta
x+1)=\ell\<a,b\>[D,(x^{3}-2\beta x^{2}+x)]=0,$$ which implies
$\ell=0$.

Furthermore, with $\ell= 0$ we have
\begin{eqnarray*}
a^{1}(x)_{0}b^{1}(x)=\Res_{x_{1}}[a^{1}(x_{1}),b^{1}(x)]=(x^{3}-2\beta
x^{2}+x)[a,b](x).
\end{eqnarray*}
Using the same reasoning we get $(3x^{2}-4\beta x+1)[a,b](x)=0$. As
$[a,b](x)\in \Hom (W,W((x)))$, it follows that $[a,b](x)=0$. With
$[\g,\g]=\g$ and $\<\cdot,\cdot\>$ non-degenerate, we have $u(x)=0$
on $W$ for every $u\in \g$. We also have $u^{1}(x)=0$ on $W$. Since
$w_{0}$ generates $W$, $W$ must be a $1$-dimensional trivial
$\hat{\g}_{e}$-module. This completes the proof.
\end{proof}

\bl{lpoly} Let $W$ be a restricted $\hat{\g}_{e}$-module of level
$(\ell,\ell_{+},\ell_{-})$ with $\ell\ne 0$. Then the vertex algebra
$V_{W}$ generated by $U_{W}$ is a $\C[x]$-submodule of $\E(W)$. \el

\begin{proof}  Let $a,b\in \g$ with $\<a,b\>\ne 0$.
We have
$$a^{1}(x)_{1}b^{1}(x)=\Res_{x_{1}}(x_{1}-x)[a^{1}(x_{1}),b^{1}(x)]
=\ell \<a,b\>(x^{3}-2\beta x^{2}+x).$$ It follows that $x^{3}-2\beta
x^{2}+x\in V_{W}$. Furthermore, we have
$$(x^{3}-2\beta x^{2}+x)_{-3}1_{W}
=\frac{1}{2}\left(\frac{d}{dx}\right)^{2}(x^{3}-2\beta
x^{2}+x)=3x-2\beta,$$ which implies $x\in V_{W}$. Noticing that
$x_{-1}u=xu$ for $u\in V_{W}$, we obtain $\C[x]V_{W}\subset V_{W}$,
proving that $V_{W}$ is a $\C[x]$-submodule of $\E(W)$.
\end{proof}

To better describe the vertex algebras generated by the generating
functions of elliptic affine Lie algebras on restricted modules, we
shall make use of certain closely related Lie algebras.

\bp{pbig} Let $\g$ be a (possibly infinite-dimensional) Lie algebra
over $\C$, equipped with a non-degenerate symmetric invariant
bilinear form $\<\cdot,\cdot\>$, and let $p(\xi)\in \C((\xi))$ (with
$\xi$ a new formal variable). Set
$$K=\C((\xi))(\g\oplus \g^{1}\oplus \C \k)\otimes
\C[t,t^{-1}],$$  a vector space over $\C$. Define
\begin{eqnarray*}
&&[f(\xi)\k\otimes t^{m}, K]=0=[K,f(\xi)\k\otimes t^{m}],\\
&&[f(\xi)a\otimes t^{m},g(\xi)b\otimes
t^{n}]=f(\xi)g(\xi)[a,b]\otimes
t^{m+n}+\<a,b\>f'(\xi)g(\xi)\k\otimes t^{m+n}\\
&&\hspace{5cm}
+m\<a,b\>f(\xi)g(\xi)\k\otimes t^{m+n-1},\\
&&[f(\xi)a\otimes t^{m}, g(\xi)b^{1}\otimes
t^{n}]=f(\xi)g(\xi)[a,b]^{1}\otimes t^{m+n},\\
&&[f(\xi)a^{1}\otimes t^{m}, g(\xi)b\otimes
t^{n}]=f(\xi)g(\xi)[a,b]^{1}\otimes t^{m+n},\\
&&[f(\xi)a^{1}\otimes t^{m}, g(\xi)b^{1}\otimes t^{n}]
=f(\xi)g(\xi)p(\xi)[a,b]\otimes t^{m+n}\\
&&\hspace{5cm}+\<a, b\>f'(\xi)g(\xi)p(\xi)\k\otimes t^{m+n}\\
&&\hspace{5cm}+\frac{1}{2}\<a, b\>f(\xi)g(\xi)p'(\xi)\k\otimes
t^{m+n}\\
&&\hspace{5cm}+m\<a, b\>f(\xi)g(\xi)p(\xi)\k\otimes t^{m+n-1}
\end{eqnarray*}
for $a,b\in \g,\; f(\xi),g(\xi)\in \C((\xi)),\; m,n\in \Z$. Let
$J_{0}$ be the subspace of $K$, spanned over $\C$ by
$$f'(\xi)\k\otimes t^{n}+nf(\xi)\k\otimes t^{n-1}$$
for $f(\xi)\in \C((\xi)),\; n\in \Z$. Then $J_{0}$ is a two-sided
ideal of the algebra $K$, and $K/J_{0}$ is a Lie algebra, which we
denote by $K(\g,p)$. Furthermore, the $\C$-linear operator
$\tilde{D}$ on $K$, defined by
$$\tilde{D}(u\otimes t^{n})=-n(u\otimes t^{n-1})$$ for $u\in
\C((\xi))(\g+\g^{1}+\C \k),\; n\in \Z$, reduces to a derivation $D$
of the Lie algebra $K(\g,p)$. \ep

\begin{proof} It is immediate that
$J_{0}$ is a two-sided ideal of $K$. We now prove that skew symmetry
holds for the quotient algebra $K/J_{0}$. We shall just consider the
two cases which are not immediate.   Let $a, b\in\g$, $f(\xi),
g(\xi)\in\C((\xi)),\ m, n\in\Z$. With $[a,b]=-[b,a]$ and
$\<a,b\>=\<b,a\>$, we have
\begin{eqnarray*}
&&[f(\xi)a\otimes t^{m}, g(\xi)b\otimes t^{n}]+[g(\xi)b\otimes t^{n},f(\xi)a\otimes t^{m}]\\
&=&\<a, b\>f'(\xi)g(\xi){\bf
k}\otimes t^{m+n}+m\<a, b\>f(\xi)g(\xi){\bf k}\otimes t^{m+n-1}\\
&&+\<a,b\>f(\xi)g'(\xi){\bf
k}\otimes t^{m+n}+n\<a, b\>f(\xi)g(\xi){\bf k}\otimes t^{m+n-1}\\
&=&\<a, b\>\left(\frac{d}{d\xi}(f(\xi)g(\xi)){\bf k}\otimes t^{m+n}
+(m+n)f(\xi)g(\xi){\bf k}\otimes t^{m+n-1}\right)\\
&\in & J_{0},
\end{eqnarray*}
\begin{eqnarray*}
&&[f(\xi)a^{1}\otimes t^{m}, g(\xi)b^{1}\otimes
t^{n}]+[g(\xi)b^{1}\otimes t^{n},f(\xi)a^{1}\otimes t^{m}]\\
&=&\<a, b\>f'(\xi)g(\xi)p(\xi){\bf k}\otimes t^{m+n} +
\<a, b\>f(\xi)g'(\xi)p(\xi){\bf k}\otimes t^{m+n}\\
&&+\frac{1}{2}\<a, b\>f(\xi)g(\xi)p'(\xi){\bf k}\otimes t^{m+n}+
\frac{1}{2}\<a, b\>f(\xi)g(\xi)p'(\xi){\bf k}\otimes t^{m+n}\\
&&+m\<a, b\>f(\xi)g(\xi)p(\xi){\bf k}\otimes t^{m+n-1} +n\<a,
b\>f(\xi)g(\xi)p(\xi){\bf k}\otimes t^{m+n-1}\\
&=&\<a,b\>\left(\frac{d}{d\xi}(f(\xi)g(\xi)p(\xi)){\bf k}\otimes
t^{m+n}+(m+n)f(\xi)g(\xi)p(\xi){\bf k}\otimes t^{m+n-1}\right)\\
 &\in& J_{0}.
\end{eqnarray*}

Next, we establish Jacobi identity. Let $f(\xi),g(\xi),h(\xi)\in
\C((\xi)),\ a,b,c\in \g,\ l,m,n\in \Z$. We have
\begin{eqnarray*}
&&[f(\xi)a\otimes t^{l}, [g(\xi)b\otimes t^{m}, h(\xi)c\otimes
t^{n}]]\\
&=&[f(\xi)a\otimes t^{l}, g(\xi)h(\xi)[b, c]\otimes t^{m+n}]\\
&=&f(\xi)g(\xi)h(\xi)[a, [b, c]]\otimes t^{l+m+n}+\<a, [b,
c]\>f'(\xi)g(\xi)h(\xi){\bf k}\otimes t^{l+m+n}\\
&&+l\<a, [b, c]\>f(\xi)g(\xi)h(\xi){\bf k}\otimes t^{l+m+n-1},
\end{eqnarray*}
\begin{eqnarray*}
&&[g(\xi)b\otimes t^{m}, [f(\xi)a\otimes t^{l}, h(\xi)c\otimes
t^{n}]]\\
&=&f(\xi)g(\xi)h(\xi)[b, [a, c]]\otimes t^{l+m+n}+\<b, [a,
c]\>f(\xi)g'(\xi)h(\xi){\bf k}\otimes t^{l+m+n}\\
&&+m\<b, [a, c]\>f(\xi)g(\xi)h(\xi){\bf k}\otimes t^{l+m+n-1},
\end{eqnarray*}
\begin{eqnarray*}
 &&[[f(\xi)a\otimes t^{l}, g(\xi)b\otimes t^{m}],
h(\xi)c\otimes
t^{n}]\\
&=&[f(\xi)g(\xi)[a, b]\otimes t^{l+m}, h(\xi)c\otimes t^{n}]\\
&=&f(\xi)g(\xi)h(\xi)[[a, b], c]\otimes t^{l+m+n}+\<[a, b],
c\>(f(\xi)g(\xi))'h(\xi){\bf k}\otimes t^{l+m+n}\\
&&+(l+m)\<[a, b], c\>f(\xi)g(\xi)h(\xi){\bf k}\otimes t^{l+m+n-1}.
\end{eqnarray*}
Then the Jacobi identity for the triple $(f(\xi)a\otimes
t^{l},g(\xi)b\otimes t^{m},h(\xi)c\otimes t^{n})$ follows from the
Jacobi identity of $\g$ and the invariance of $\<\cdot, \cdot\>$.

We also have
\begin{eqnarray*}
[f(\xi)a\otimes t^{l}, [g(\xi)b\otimes t^{m}, h(\xi)c^{1}\otimes
t^{n}]]
&=&[f(\xi)a\otimes t^{l}, g(\xi)h(\xi)[b, c]^{1}\otimes t^{m+n}]\\
&=&f(\xi)g(\xi)h(\xi)[a, [b, c]]^{1}\otimes t^{l+m+n},
\end{eqnarray*}
\begin{eqnarray*}
[g(\xi)b\otimes t^{m}, [f(\xi)a\otimes t^{l}, h(\xi)c^{1}\otimes
t^{n}]] =f(\xi)g(\xi)h(\xi)[b, [a, c]]^{1}\otimes t^{l+m+n},
\end{eqnarray*}
\begin{eqnarray*}
[[f(\xi)a\otimes t^{l}, g(\xi)b\otimes t^{m}], h(\xi)c^{1}\otimes
t^{n}]
&=&[f(\xi)g(\xi)[a, b]\otimes t^{l+m}, h(\xi)c^{1}\otimes t^{n}]\\
&=&f(\xi)g(\xi)h(\xi)[[a, b], c]^{1}\otimes t^{l+m+n}.
\end{eqnarray*}
These imply the Jacobi identity for the indicated triple.

We have
\begin{eqnarray*}
&&[f(\xi)a\otimes t^{l}, [g(\xi)b^{1}\otimes t^{m},
h(\xi)c^{1}\otimes
t^{n}]]\\
&=&[f(\xi)a\otimes t^{l}, g(\xi)h(\xi)p(\xi)[b, c]\otimes t^{m+n}]\\
&=&f(\xi)g(\xi)h(\xi)p(\xi)[a, [b, c]]\otimes t^{l+m+n}+\<a, [b,
c]\>f'(\xi)g(\xi)h(\xi)p(\xi){\bf k}\otimes t^{l+m+n}\\
&&+l\<a, [b, c]\>f(\xi)g(\xi)h(\xi)p(\xi){\bf k}\otimes t^{l+m+n-1},
\end{eqnarray*}
\begin{eqnarray*}
&&[g(\xi)b^{1}\otimes t^{m}, [f(\xi)a\otimes t^{l},
h(\xi)c^{1}\otimes
t^{n}]]\\
&=&[g(\xi)b^{1}\otimes t^{m}, f(\xi)h(\xi)[a, c]^{1}\otimes t^{l+m}]\\
&=&f(\xi)g(\xi)h(\xi)p(\xi)[b, [a, c]]\otimes t^{l+m+n}+\<b, [a,
c]\>f(\xi)g'(\xi)h(\xi)p(\xi){\bf k}\otimes t^{l+m+n}\\
&&+\frac{1}{2}\<b, [a, c]\>f(\xi)g(\xi)h(\xi)p'(\xi){\bf k}\otimes
t^{l+m+n}\\
&&+m\<b, [a, c]\>f(\xi)g(\xi)h(\xi){\bf k}\otimes t^{l+m+n-1},
\end{eqnarray*}
\begin{eqnarray*}
 &&[[f(\xi)a\otimes t^{l}, g(\xi)b^{1}\otimes t^{m}],
h(\xi)c^{1}\otimes
t^{n}]\\
&=&[f(\xi)g(\xi)[a, b]^{1}\otimes t^{l+m}, h(\xi)c^{1}\otimes t^{n}]\\
&=&f(\xi)g(\xi)h(\xi)p(\xi)[[a, b], c]\otimes t^{l+m+n}\\
&&+\<[a, b],
c\>(f'(\xi)g(\xi)+f(\xi)g'(\xi))h(\xi)p(\xi){\bf k}\otimes t^{l+m+n}\\
&&+\frac{1}{2}\<[a, b], c\>f(\xi)g(\xi)h(\xi)p'(\xi){\bf k}\otimes
t^{l+m+n}\\
&&+(l+m)\<[a, b], c\>f(\xi)g(\xi)h(\xi)p(\xi){\bf k}\otimes
t^{l+m+n-1}.
\end{eqnarray*}
Using the Jacobi identity of $\g$ and the invariance of $\<\cdot,
\cdot\>$ we obtain the Jacobi identity for the indicated triple.

We have
\begin{eqnarray*}
&&[f(\xi)a^{1}\otimes t^{l}, [g(\xi)b^{1}\otimes t^{m},
h(\xi)c^{1}\otimes
t^{n}]]\\
&=&[f(\xi)a^{1}\otimes t^{l}, g(\xi)h(\xi)p(\xi)[b, c]\otimes t^{m+n}]\\
&=&f(\xi)g(\xi)h(\xi)p(\xi)[a, [b, c]]^{1}\otimes t^{l+m+n},
\end{eqnarray*}
\begin{eqnarray*}
[g(\xi)b^{1}\otimes t^{m}, [f(\xi)a^{1}\otimes t^{l},
h(\xi)c^{1}\otimes t^{n}]]=f(\xi)g(\xi)h(\xi)p(\xi)[b, [a,
c]]^{1}\otimes t^{l+m+n},
\end{eqnarray*}
\begin{eqnarray*}
[[f(\xi)a^{1}\otimes t^{l}, g(\xi)b^{1}\otimes t^{m}],
h(\xi)c^{1}\otimes t^{n}]=f(\xi)g(\xi)h(\xi)p(\xi)[[a, b],
c]^{1}\otimes t^{l+m+n}.
\end{eqnarray*}
Then the Jacobi identity for the indicated triple follows.
Furthermore, the other cases follow from these and skew symmetry.

As for the last assertion, we first show that $\tilde{D}$ is a
derivation of algebra $K$, by checking the two nontrivial cases. We
have
\begin{eqnarray*}
&&-m[f(\xi)a\otimes t^{m-1}, g(\xi)b\otimes t^{n}]-n[f(\xi)a\otimes
t^{m}, g(\xi)b\otimes t^{n-1}]\\
&=&-mf(\xi)g(\xi)[a,b]\otimes
t^{m+n-1}-m\<a,b\>f'(\xi)g(\xi)\k\otimes t^{m+n-1}\\
&&\hspace{5cm} -m(m-1)\<a,b\>f(\xi)g(\xi)\k\otimes t^{m+n-2}\\
&&-nf(\xi)g(\xi)[a,b]\otimes
t^{m+n-1}-n\<a,b\>f'(\xi)g(\xi)\k\otimes t^{m+n-1}\\
&&\hspace{5cm} -mn\<a,b\>f(\xi)g(\xi)\k\otimes t^{m+n-2}\\
&=&-(m+n)f(\xi)g(\xi)[a,b]\otimes t^{m+n-1}
-\<a,b\>(m+n)f'(\xi)g(\xi)\k\otimes t^{m+n-1}\\
&&-m(m+n-1)\<a,b\>f(\xi)g(\xi)\k\otimes t^{m+n-2},
\end{eqnarray*}
\begin{eqnarray*}
&&-m[f(\xi)a^{1}\otimes t^{m-1}, g(\xi)b^{1}\otimes t^{n}]
-n[f(\xi)a^{1}\otimes t^{m}, g(\xi)b^{1}\otimes t^{n-1}]\\
&=&-mf(\xi)g(\xi)p(\xi)[a,b]\otimes t^{m+n-1}
-m\<a, b\>f'(\xi)g(\xi)p(\xi)\k\otimes t^{m+n-1}\\
&&-\frac{1}{2}\<a, b\>mf(\xi)g(\xi)p'(\xi)\k\otimes
t^{m+n-1}-m(m-1)\<a, b\>f(\xi)g(\xi)p(\xi)\k\otimes t^{m+n-2}\\
&&-nf(\xi)g(\xi)p(\xi)[a,b]\otimes t^{m+n-1}
-n\<a, b\>f'(\xi)g(\xi)p(\xi)\k\otimes t^{m+n-1}\\
&&-\frac{1}{2}\<a, b\>nf(\xi)g(\xi)p'(\xi)\k\otimes t^{m+n-1}-mn\<a,
b\>f(\xi)g(\xi)p(\xi)\k\otimes t^{m+n-2}\\
&=&-(m+n)f(\xi)g(\xi)p(\xi)[a,b]\otimes t^{m+n-1}
-(m+n)\<a, b\>f'(\xi)g(\xi)p(\xi)\k\otimes t^{m+n-1}\\
&&-\frac{1}{2}\<a, b\>(m+n)f(\xi)g(\xi)p'(\xi)\k\otimes t^{m+n-1}\\
&&-m(m+n-1)\<a, b\>f(\xi)g(\xi)p(\xi)\k\otimes t^{m+n-2}.
\end{eqnarray*}
It is readily seen that $\tilde{D}$ preserves the ideal $J_{0}$, so
that $\tilde{D}$ reduces to a derivation $D$ of $K(\g,p)$.
\end{proof}

Note that though $K(\g,p)$ is naturally a vector space over the
field $\C((\xi))$, $K(\g,p)$ is {\em not} a Lie algebra over
$\C((\xi))$, as the Lie bracket is not $\C((\xi))$-bilinear.

We see from the construction that as a vector space over $\C$,
\begin{eqnarray}
K(\g,p)=(\C((\xi))(\g\oplus \g^{1})\otimes
\C[t,t^{-1}])\oplus(R/dR),
\end{eqnarray}
where $R=\C((\xi))\k\otimes \C[t,t^{-1}]$ and
$d=\frac{d}{d\xi}\otimes 1+1\otimes \frac{d}{dt}$. For $u\in
\C((\xi))(\g+\g^{1}+\C \k), \; n\in \Z$, denote by $u(n)$ the image
of $u\otimes t^{n}$ in $K(\g,p)$. Since
$$d(\k\otimes t^{n})=n(\k\otimes t^{n-1})\ \ \ \mbox{ for }n\in
\Z,$$ we have
\begin{eqnarray}
{\bf k}(n)=0\ \ \ \mbox{ for }n\ne -1.
\end{eqnarray}

For $u\in \C((\xi))(\g+\g^{1}+\C \k)$, form the generating function
\begin{eqnarray}
Y_{t}(u,x)=\sum_{n\in \Z}u(n)x^{-n-1}\in K(\g,p)[[x,x^{-1}]].
\end{eqnarray}
The Lie bracket relations become
\begin{eqnarray}\label{ekgp-relations}
&&[Y_{t}(f\k,x_{1}), K(\g,p)]=0,\nonumber\\
&&[Y_{t}(fa,x_{1}),Y_{t}(gb,x_{2})]
=Y_{t}(fg[a,b],x_{2})x_{2}^{-1}\delta\left(\frac{x_{1}}{x_{2}}\right)\nonumber\\
&&\ \ \ \
+\<a,b\>Y_{t}(f'g\k,x_{2})x_{2}^{-1}\delta\left(\frac{x_{1}}{x_{2}}\right)
+\<a,b\>Y_{t}(fg\k,x_{2})\frac{\partial}{\partial
x_{2}}x_{2}^{-1}\delta\left(\frac{x_{1}}{x_{2}}\right),\ \ \nonumber\\
&&[Y_{t}(fa,x_{1}),Y_{t}(gb^{1},x_{2})]
=Y_{t}(fg[a,b]^{1},x_{2})x_{2}^{-1}\delta\left(\frac{x_{1}}{x_{2}}\right),\nonumber\\
&&[Y_{t}(fa^{1},x_{1}),Y_{t}(gb,x_{2})]
=Y_{t}(fg[a,b]^{1},x_{2})x_{2}^{-1}\delta\left(\frac{x_{1}}{x_{2}}\right),\nonumber\\
&&[Y_{t}(fa^{1},x_{1}),Y_{t}(gb^{1},x_{2})]=Y_{t}(fgp[a,b],x_{2})
x_{2}^{-1}\delta\left(\frac{x_{1}}{x_{2}}\right)\nonumber\\
&&\ \ \ \ +\frac{1}{2}\<a,b\>Y_{t}((2f'gp+fgp')\k,
x_{2})x_{2}^{-1}\delta\left(\frac{x_{1}}{x_{2}}\right)\nonumber\\
&&\ \ \ \ +\<a,b\>Y_{t}(fgp\k,x_{2})\frac{\partial}{\partial
x_{2}}x_{2}^{-1}\delta\left(\frac{x_{1}}{x_{2}}\right).
\end{eqnarray}

{}From now on, we assume that $p(x)\ (\in \C[x])$ is a polynomial.
Set
\begin{eqnarray}
K^{0}(\g,p)=(\C[\xi](\g\oplus \g^{1})\otimes
\C[t,t^{-1}])\oplus(R/dR)\subset K,
\end{eqnarray}
where $R$ and $d$ are given as before. It is readily seen that
$K^{0}(\g,p)$ is a Lie subalgebra which is stable under the
derivation $D$.

We now construct a family of Lie algebras, generalizing the elliptic
affine Lie algebra $\hat{\g}_{e}$ (with $\k_{\pm}=0$).

\bp{phatgp} Let $\g$ be a (possibly infinite-dimensional) Lie
algebra over $\C$, equipped with a non-degenerate symmetric
invariant bilinear form $\<\cdot,\cdot\>$, and let $p(x)\in \C[x]$.
Set
\begin{eqnarray}
\hat{\g}_{p}=(\g\oplus \g^{1})\otimes \C[t,t^{-1}]\oplus \C \k,
\end{eqnarray}
a vector space over $\C$. Define a bilinear operation on
$\hat{\g}_{p}$ by
\begin{eqnarray}\label{ep-elliptic}
&&[\k, \hat{\g}_{p}]=0=[\hat{\g}_{p},\k],\nonumber\\
&&[a(x_{1}),b(x_{2})]
=[a,b](x_{2})x_{2}^{-1}\delta\left(\frac{x_{1}}{x_{2}}\right)
+\<a,b\>\k\frac{\partial}{\partial
x_{2}}x_{2}^{-1}\delta\left(\frac{x_{1}}{x_{2}}\right),\ \ \nonumber\\
&&[a(x_{1}),b^{1}(x_{2})]
=[a,b]^{1}(x_{2})x_{2}^{-1}\delta\left(\frac{x_{1}}{x_{2}}\right),\nonumber\\
&&[a^{1}(x_{1}),b(x_{2})]
=[a,b]^{1}(x_{2})x_{2}^{-1}\delta\left(\frac{x_{1}}{x_{2}}\right),\nonumber\\
&&[a^{1}(x_{1}),b^{1}(x_{2})]=p(x_{2})[a,b](x_{2})
x_{2}^{-1}\delta\left(\frac{x_{1}}{x_{2}}\right)\nonumber\\
&&\hspace{1cm} \ \ +\frac{1}{2}\<a,b\>p'(x_{2})\k
x_{2}^{-1}\delta\left(\frac{x_{1}}{x_{2}}\right)
+\<a,b\>p(x_{2})\k\frac{\partial}{\partial
x_{2}}x_{2}^{-1}\delta\left(\frac{x_{1}}{x_{2}}\right)\!.\ \ \ \ \ \
\
\end{eqnarray}
Then $\hat{\g}_{p}$ is a Lie algebra over $\C$. \ep

\begin{proof} Recall the Lie algebra
$$K^{0}(\g,p)= (\C[\xi](\g\oplus \g^{1})\otimes
\C[t,t^{-1}])\oplus(R/dR)\subset K. $$ Let $J$ be the $\C$-span of
the coefficients of
$$(f(\xi)a)(x)-f(x)a(x),\ \ \ (f(\xi)a^{1})(x)-f(x)a^{1}(x), \ \ \
(g(\xi)\k)(x)-g(x)\k$$ for $a\in \g,\; f(\xi)\in \C[\xi],\;
g(\xi)\in \C((\xi))$. By using (\ref{ekgp-relations}), it is
straightforward to show that $J$ is a left ideal of $K^{0}(\g,p)$.
We see that the underlying vector space of the quotient Lie algebra
$K^{0}(\g,p)/J$ is isomorphic to $\hat{\g}_{p}$. Then it follows
immediately that the defined nonassociative algebra $\hat{\g}_{p}$
is a Lie algebra.
\end{proof}

\br{relliptic-lie} {\em Note that the Lie algebra $\hat{\g}_{p}$
with $p(x)=x^{3}-2\beta x^{2}+x$ is isomorphic to the quotient
algebra $\hat{\g}_{e}/(\C\k_{+}+\C \k_{-})$ of the elliptic affine
Lie algebra.} \er

We next construct another  family of Lie algebras.

\bp{pcheckgp} Let $\g$ be a (possibly infinite-dimensional) Lie
algebra over $\C$, equipped with a non-degenerate symmetric
invariant bilinear form $\<\cdot,\cdot\>$, and let $p(x)\in \C[x]$.
Set
\begin{eqnarray}
\check{\g}_{p}=\C((z))\otimes (\g\oplus \g^{1})\otimes
\C[t,t^{-1}]\oplus \C((z)) \k,
\end{eqnarray}
a vector space over $\C((z))$. For $a\in \g$, set
 $$a(x)=\sum_{n\in \Z} (a\otimes t^{n})x^{-n-1},\ \ \
a^{1}(x)=\sum_{n\in \Z} (a^{1}\otimes t^{n})x^{-n-1}.$$ Define a
$\C((z))$-bilinear operation on $\check{\g}_{p}$ by
\begin{eqnarray}\label{ekgp-relations-final}
&&[\k, \check{\g}_{p}]=0=[\check{\g}_{p},\k],\nonumber\\
&&[a(x_{1}),b(x_{2})]
=[a,b](x_{2})x_{2}^{-1}\delta\left(\frac{x_{1}}{x_{2}}\right)
+\<a,b\>\k\frac{\partial}{\partial
x_{2}}x_{2}^{-1}\delta\left(\frac{x_{1}}{x_{2}}\right),\ \ \nonumber\\
&&[a(x_{1}),b^{1}(x_{2})]
=[a,b]^{1}(x_{2})x_{2}^{-1}\delta\left(\frac{x_{1}}{x_{2}}\right),\nonumber\\
&&[a^{1}(x_{1}),b(x_{2})]
=[a,b]^{1}(x_{2})x_{2}^{-1}\delta\left(\frac{x_{1}}{x_{2}}\right),\nonumber\\
&&[a^{1}(x_{1}),b^{1}(x_{2})]=p(z+x_{2})[a,b](x_{2})
x_{2}^{-1}\delta\left(\frac{x_{1}}{x_{2}}\right)\nonumber\\
&&\ \ \ +\frac{1}{2}\<a,b\>p'(z+x_{2})\k
x_{2}^{-1}\delta\left(\frac{x_{1}}{x_{2}}\right)
+\<a,b\>p(z+x_{2})\k\frac{\partial}{\partial
x_{2}}x_{2}^{-1}\delta\left(\frac{x_{1}}{x_{2}}\right)\!. \ \ \ \ \
\ \
\end{eqnarray}
Then $\check{\g}_{p}$ is a Lie algebra over $\C((z))$. Furthermore,
the map $\check{D}$, defined by
\begin{eqnarray}
\check{D}(f(z)\k)=f'(z)\k,\ \ \ \check{D}\left(f(z)u\otimes
t^{n}\right)=f'(z)u\otimes t^{n}-nf(z)u\otimes t^{n-1}
\end{eqnarray}
for $f(z)\in \C((z)),\; u\in \g+\g^{1},\; n\in \Z$, is a $\C$-linear
derivation of $\check{\g}_{p}$ viewed as a Lie algebra over $\C$.\ep

\begin{proof} Just as in the proof of Proposition \ref{phatgp}
we shall make use of the Lie algebra $K^{0}(\g,p)$. We extend
$K^{0}(\g,p)$ to a Lie algebra $\C((z))\otimes_{\C} K^{0}(\g,p)$
over the filed $\C((z))$. Let $J_{1}$ be the subspace of
$\C((z))\otimes K^{0}(\g,p)$, spanned over $\C((z))$ by the
coefficients of
$$(f(\xi)u)(x)-f(z+x)u(x), \ \ \
(g(\xi)\k)(x)-g(z+x)\k$$ for $a\in \g+\g^{1},\; f(\xi)\in \C[\xi],\;
g(\xi)\in \C((\xi))$. By using (\ref{ekgp-relations}), it is
straightforward to show that $J_{1}$ is a left ideal of
$\C((z))\otimes_{\C} K^{0}(\g,p)$. One sees that the underlying
vector space of the quotient Lie algebra $(\C((z))\otimes_{\C}
K^{0}(\g,p))/J_{1}$ is isomorphic to $\check{\g}_{p}$. Then the
first assertion follows immediately.

As for the furthermore assertion, recall that $D$ is a derivation of
$K^{0}(\g,p)$. Then $\frac{d}{dz}\otimes 1+1\otimes D$ is a
$\C$-linear derivation of $\C((z))\otimes_{\C} K^{0}(\g,p)$ viewed
as a Lie algebra over $\C$. We have
\begin{eqnarray*}
&&\left(\frac{d}{dz}\otimes 1+1\otimes D\right)\left(1\otimes
f(\xi)u\otimes t^{n}-\sum_{j\ge 0}\frac{1}{j!}f^{(j)}(z)\otimes
u\otimes t^{n+j}\right)\\
&=&-n(1\otimes f(\xi)u\otimes t^{n-1}) +\sum_{j\ge
0}(n+j)\frac{1}{j!}f^{(j)}(z)\otimes
u\otimes t^{n+j-1}\\
&&-\sum_{j\ge 0}\frac{1}{j!}f^{(j+1)}(z)\otimes u\otimes
t^{n+j}\\
&=&-n\left(1\otimes f(\xi)u\otimes t^{n-1}-\sum_{j\ge
0}\frac{1}{j!}f^{(j)}(z)\otimes u\otimes t^{n+j-1}\right)
\end{eqnarray*}
for $u\in \g+\g^{1},\ f(\xi)\in \C[\xi],\; n\in \Z,$ or for $u\in
\C\k, f(\xi)\in \C((\xi)),\; n\in \Z$. It follows that $J_{1}$ is
stable under $\frac{d}{dz}\otimes 1+1\otimes D$. Then
$\frac{d}{dz}\otimes 1+1\otimes D$ gives rise to a derivation of
$\check{\g}_{p}$ viewed as a Lie algebra  over $\C$.
\end{proof}

\bd{delliptic-tva} {\em  We define $\check{\g}_{e}$ to be the Lie
algebra $\check{\g}_{p}$ with $p(x)=x^{3}-2\beta x^{2}+x$.} \ed

Lie algebra $\check{\g}_{p}$ as a vector space is naturally
$\Z$-graded
$$\check{\g}_{p}=\coprod_{n\in \Z}
\left(\C((z))(\g+\g^{1})\otimes
t^{n}+\delta_{n,0}\C((z))\k\right),$$ but this does {\em not} make
$\check{\g}_{p}$ a $\Z$-graded Lie algebra. Nevertheless, we can
make $\check{\g}_{p}$ a $\Z$-filtered Lie algebra. The following is
straightforward from the defining relations of $\check{\g}_{p}$ and
{}from the assumption that $p(x)$ is a polynomial:

\bl{lfilter} Let $n\in \Z$. Set
\begin{eqnarray}
\check{\g}_{p}[n]=
\begin{cases}\C((z))(\g+\g^{1})\otimes
t^{n}\C[t] & \mbox{ for } n\ge 1,\\
\C((z))(\g+\g^{1})\otimes t^{n}\C[t]\oplus \C((z))\k & \mbox{ for
}n\le 0.
\end{cases}
\end{eqnarray}
Then $\{\check{\g}_{p}[n]\}_{n\in \Z}$ is a decreasing filtration of
$\check{\g}_{p}$, satisfying that
\begin{eqnarray}
&&\cap_{n\in \Z}\check{\g}_{p}[n]=0,\\
&&[\check{\g}_{p}[m],\check{\g}_{p}[n]]\subset \check{\g}_{p}[m+n] \
\ \ \mbox{ for }m,n\in \Z.
\end{eqnarray}
\el

\bd{dvacuum-checkmodule} {\em We say that a $\check{\g}_{p}$-module
$W$ is of {\em level} $\ell\in \C$ if $\k$ acts on $W$ as scalar
$\ell$. We define a {\em vacuum $\check{\g}_{p}$-module} to be a
$\check{\g}_{p}$-module $W$ equipped with a vector $w_{0}\in W$ and
a $\C$-linear operator $D$ on $W$ such that
$W=U(\check{\g}_{p})w_{0}$, $Dw_{0}=0$,
$$u(n)w_{0}=0\ \ \ \mbox{ for }u\in \g+\g^{1},\; n\ge 0,$$ and such
that
$$[D,u(x)]=\frac{d}{dx}u(x),\ \ \ \
[D,f(z)]=f'(z)$$ for $u\in \g+\g^{1},\;f(z)\in \C((z))$.} \ed

\bl{lvacuum-restricted}
Let $W$ be a $\check{\g}_{p}$-module with a
vector $w_{0}$ satisfying that
$$W=U(\check{\g}_{p})w_{0}\  \mbox{ and } \ u(n)w_{0}=0\ \ \ \mbox{
for }u\in \g+\g^{1},\; n\ge 0.$$ Then $W$ is restricted. In
particular,  any vacuum $\check{\g}_{p}$-module is a restricted
module.\el

\begin{proof} We need to prove that for any $w\in W$,
$\check{\g}_{p}[n]w=0$ for $n$ sufficiently large. For $k\ge 0$, let
$W[k]$ be the span of the subspaces
$$\check{\g}_{p}[n_{1}]\cdots \check{\g}_{p}[n_{r}]w_{0}$$
for $0\le r\le k$ with $n_{1},\dots, n_{r}\in \Z$. Then
$W=\cup_{k\ge 0}W[k]$.
 {}From definition we have $W[0]=\C((z))w_{0}$, so we have
$\check{\g}_{p}[n]W[0]=0$ for $n\ge 1$. It follows from Lemma
\ref{lfilter} and induction on $k$ that for any $k\ge 0$ and for any
$w\in W[k]$, $\check{\g}_{p}[n]w=0$ for $n$ sufficiently large. As
$W=\cup_{k\ge 0}W[k]$, it follows that $W$ is restricted.
\end{proof}

\bp{pvacuum-module} Let $\ell$ be a complex number and let $W$ be a
restricted $\hat{\g}_{p}$-module of level $\ell$. Set
$V_{W}=\C((x))\<U_{W}\>$, where
$$U_{W}=\{ a(x),a^{1}(x)\;|\; a\in \g\}\subset \E(W).$$
Then $V_{W}$ with vector $1_{W}\in V_{W}$ and operator $D=\D$ (the
$\D$-operator of $V_{W}$) is a vacuum $\check{\g}_{p}$-module of
level $\ell$ with $u(x_{0})$ acting as $Y_{\E}(u(x),x_{0})$ for
$u\in \g+\g^{1}$ and with $f(z)\in \C((z))$ acting as $f(x)$. \ep

\begin{proof} Note that $W$ is a faithful module for $V_{W}$
viewed as a vertex algebra over $\C$ with
$Y_{W}(\psi(x),x_{0})=\psi(x_{0})$ for $\psi(x)\in V_{W}$. Let $a,
b\in \g$. We have
\begin{eqnarray*}
&&[Y_{W}(a(x),x_{1}), Y_{W}(b(x),x_{2})]\\
&=&Y_{W}([a,
b](x),x_{2})x_{2}^{-1}\delta\left(\frac{x_{1}}{x_{2}}\right)
+\ell\<a,b\>\frac{\partial}{\partial
x_{2}}x_{2}^{-1}\delta\left(\frac{x_{1}}{x_{2}}\right),
\end{eqnarray*}
\begin{eqnarray*}
[Y_{W}(a(x),x_{1}),Y_{W}(b^{1}(x),x_{2})]
=Y_{W}([a,b]^{1}(x),x_{2})x_{2}^{-1}\delta\left(\frac{x_{1}}{x_{2}}\right),
\end{eqnarray*}
\begin{eqnarray*}
&&[Y_{W}(a^{1}(x),x_{1}),Y_{W}(b^{1}(x),x_{2})]\
\left(=[a^{1}(x_{1}), b^{1}(x_{2})]\right)\\
&=&p(x_{2})[a,b](x_{2})
x_{2}^{-1}\delta\left(\frac{x_{1}}{x_{2}}\right)
+\ell\<a,b\>p(x_{2})\frac{\partial}{\partial
x_{2}}x_{2}^{-1}\delta\left(\frac{x_{1}}{x_{2}}\right)\\
&&\hspace{0.5cm}+\frac{1}{2}\ell\<a,b\>
p'(x_{2})x_{2}^{-1}\delta\left(\frac{x_{1}}{x_{2}}\right)\\
&=&Y_{W}(p(x)[a,b](x),x_{2})x_{2}^{-1}\delta\left(\frac{x_{1}}{x_{2}}\right)
+\ell\<a,b\>Y_{W}(p(x)1_{W},x_{2})\frac{\partial}{\partial
x_{2}}x_{2}^{-1}\delta\left(\frac{x_{1}}{x_{2}}\right)\\
&&\hspace{0.5cm}+\frac{1}{2}\ell\<a,b\>Y_{W}(p'(x)1_{W},x_{2})
x_{2}^{-1}\delta\left(\frac{x_{1}}{x_{2}}\right).
\end{eqnarray*}
In view of Proposition \ref{pcommutator}, we have
\begin{eqnarray*}
&&[Y_{\E}(a(x),x_{1}), Y_{\E}(b(x),x_{2})]\\
&=&Y_{\E}([a,b](x),x_{2})x_{2}^{-1}\delta\left(\frac{x_{1}}{x_{2}}\right)
+\ell \<a,b\>\frac{\partial}{\partial
x_{2}}x_{2}^{-1}\delta\left(\frac{x_{1}}{x_{2}}\right),\\
&&[Y_{\E}(a(x),x_{1}),Y_{\E}(b^{1}(x),x_{2})]
=Y_{\E}([a,b]^{1}(x),x_{2})x_{2}^{-1}\delta\left(\frac{x_{1}}{x_{2}}\right),\\
&&[Y_{\E}(a^{1}(x),x_{1}),Y_{\E}(b^{1}(x),x_{2})]\\
&=&Y_{\E}(p(x)[a,b](x),x_{2})x_{2}^{-1}\delta\left(\frac{x_{1}}{x_{2}}\right)+\ell
\<a,b\>Y_{\E}(p(x)1_{W},x_{2})\frac{\partial}{\partial
x_{2}}x_{2}^{-1}\delta\left(\frac{x_{1}}{x_{2}}\right)\\
&&\hspace{0.5cm}+\frac{1}{2}\ell \<a,b\>Y_{\E}(
p'(x)1_{W},x_{2})x_{2}^{-1}\delta\left(\frac{x_{1}}{x_{2}}\right)\\
&=&p(x+x_{2})Y_{\E}([a,b](x),x_{2})x_{2}^{-1}\delta\left(\frac{x_{1}}{x_{2}}\right)+\ell
\<a,b\>p(x+x_{2})\frac{\partial}{\partial
x_{2}}x_{2}^{-1}\delta\left(\frac{x_{1}}{x_{2}}\right)\\
&&\hspace{0.5cm}+\frac{1}{2}\ell \<a,b\>
p'(x+x_{2})x_{2}^{-1}\delta\left(\frac{x_{1}}{x_{2}}\right).
\end{eqnarray*}
Thus $V_{W}$ is a restricted $\check{\g}_{p}$-module of level $\ell$
with $u(x_{0})$ acting as $Y_{\E}(u(x),x_{0})$ for $u\in \g+\g^{1}$
and with $f(z)\in \C((z))$ acting as $f(x)$. From the definition of
$V_{W}$, $V_{W}$ is generated over $\C((x))$ from $1_{W}$ by
operators  $a(x)_{n},\; a^{1}(x)_{n},\; f(x)$ for $u\in \g,\; n\in
\Z$. It follows that $V_{W}=U(\check{\g}_{p})1_{W}$. For the vertex
algebra $V_{W}$ with the $\D$-operator $\D$, we have
$$[\D, Y_{\E}(\psi(x),x_{0})]=\frac{d}{dx_{0}}Y_{\E}(\psi(x),x_{0})
\ \ \ \mbox{ for }\psi(x)\in V_{W}.$$ In particular, we have
\begin{eqnarray*}
[\D, Y_{\E}(u(x),x_{0})]=\frac{d}{dx_{0}}Y_{\E}(u(x),x_{0}) \ \ \
\mbox{ for }u\in \g+\g^{1}.
\end{eqnarray*}
With $\D=\frac{d}{dx}$, we also have
\begin{eqnarray*}
[\D,f(x)]=f'(x)\ \ \ \mbox{ for } f(x)\in \C((x)).
\end{eqnarray*}
Therefore, $V_{W}$ is a vacuum $\check{\g}_{p}$-module of level
$\ell$.
\end{proof}

Set
\begin{eqnarray}
\check{B}=\C((z))\otimes (\g\oplus \g^{1})\otimes \C[t]\oplus
\C((z)) \k\subset \check{\g}_{p}.
\end{eqnarray}
It is readily seen that $\check{B}$ is a subalgebra which is a
direct sum of the subalgebras $\C((z))\otimes (\g\oplus
\g^{1})\otimes \C[t]$ and $\C((z)) \k$.

Let $\ell\in \C$. We make $\C((z))$ a $\check{B}$-module by letting
$\k$ act as scalar $\ell$ and letting $\C((z))\otimes (\g\oplus
\g^{1})\otimes \C[t]$ act trivially. Denote this $\check{B}$-module
by $\C((z))_{\ell}$. Then form the induced module
\begin{eqnarray}
V_{\check{\g}_{p}}(\ell,0)=U(\check{\g}_{p})\otimes_{U(\check{B})}\C((z))_{\ell},
\end{eqnarray}
a $\check{\g}_{p}$-module. Set
$${\bf 1}=1\otimes 1\in V_{\check{\g}_{p}}(\ell,0).$$
In view of the P-B-W theorem, we can and we shall identify
$\C((z))\otimes (\g\oplus \g^{1})$ as a subspace of
$V_{\check{\g}_{p}}(\ell,0)$ through the $\C((z))$-linear map
$$f(z)u\mapsto f(z)u(-1){\bf 1}.$$
We also have $\C((z)){\bf 1}\subset V_{\check{\g}_{p}}(\ell,0)$. It
is clear that $V_{\check{\g}_{p}}(\ell,0)$ is a vacuum
$\check{\g}_{p}$-module which is universal in the obvious sense.

\bt{tvtalgebra-elliptic} Let $\ell$ be any complex number. There
exists a vertex $\C((z))$-algebra structure on the
$\check{\g}_{p}$-module $V_{\check{\g}_{p}}(\ell,0)$, which is
uniquely determined by the conditions that ${\bf 1}$ is the vacuum
vector and that
\begin{eqnarray}
Y(u,x)=u(x)\ \ \ \mbox{ for }u\in \g+\g^{1}.
\end{eqnarray}
\et

\begin{proof} It is clear that $V_{\check{\g}_{p}}(\ell,0)$ as a
$\C((z))$-vector space  is generated from ${\bf 1}$ by operators
$u(n)$ for $u\in \g+\g^{1},\; n\in \Z$. It follows that
$V_{\check{\g}_{p}}(\ell,0)$ as a $\C$-vector space is generated
from ${\bf 1}$ by operators $u(n)$ for $u\in \g+\g^{1},\; n\in \Z$
and by the left multiplication of $f(z)$ for $f(z)\in \C((z))$. Then
the structure of a vertex algebra over $\C$ such that
$$Y(u,x)=u(x),\ \ \ Y(f(z){\bf 1},x)=f(z+x)\ \ \ \mbox{ for }u\in
\g+\g^{1},\; f(z)\in \C((z))$$ is unique. Consequently, the vertex
$\C((z))$-algebra structure with the required properties is unique.

For the existence, we first show that there is a vertex algebra
structure over $\C$. Clearly, the $\C$-linear derivation $\check{D}$
of $\check{\g}_{p}$ preserves $\check{B}$. Then $\check{D}$ gives
rise to a $\C$-linear operator $\D$ on $V_{\check{\g}_{p}}(\ell,0)$,
satisfying the condition that $\D{\bf 1}=0$,
\begin{eqnarray*}
&&[\D,u(x)]=\check{D}(u(x))=\frac{d}{dx}u(x)\ \ \ \mbox{ for }u\in \g+\g^{1},\\
&&[\D,f(z)]=\check{D}(f(z))=f'(z)\ \ \ \mbox{ for }f(z)\in \C((z)).
\end{eqnarray*}
{}From the construction, we have
$V_{\check{\g}_{p}}(\ell,0)=U(\check{\g}_{p}){\bf 1}$ and
$$u(x){\bf 1}\in V_{\check{\g}_{p}}(\ell,0)[[x]]
\ \mbox{ and }\lim_{x\rightarrow 0}u(x){\bf 1}=u(-1){\bf 1}=u \ \ \
\mbox{ for }u\in \g+\g^{1}.$$

Furthermore, by lemma \ref{lvacuum-restricted},
$V_{\check{\g}_{p}}(\ell,0)$ is a restricted
$\check{\g}_{p}$-module, and then it follows from the commutation
relations of $\check{\g}_{p}$ that
$$\{ f(z+x)a(x), \;
f(z+x)a^{1}(x),\; f(z+x)\; |\; a\in \g,\; f(z)\in \C((z))\}$$ is a
local subset of $\E(V_{\check{\g}_{p}}(\ell,0))$. Also, for $u\in
\g+\g^{1},\; f(z)\in \C((z))$, we have
\begin{eqnarray*}
&&f(z+x)u(x){\bf 1}\in V_{\check{\g}_{p}}(\ell,0)[[x]]\ \mbox{ and
}\
\lim_{x\rightarrow 0}f(z+x)u(x){\bf 1}=f(z)u,\\
&&f(z+x){\bf 1}\in V_{\check{\g}_{p}}(\ell,0)[[x]]\ \mbox{ and }\
\lim_{x\rightarrow 0}f(z+x){\bf 1}=f(z){\bf 1},
\end{eqnarray*}
and
\begin{eqnarray*}
&&[\D,f(z+x)u(x)]=\frac{d}{dx}\left(f(z+x)u(x)\right),\\
&&[\D,f(z+x)]=f'(z+x)=\frac{d}{dx}f(z+x).
\end{eqnarray*}
 Now, it follows {}from a
theorem of Frenkel-Kac-Radul-Wang \cite{fkrw} and Meurman-Primc
\cite{mp} that $V_{\check{\g}_{p}}(\ell,0)$ has a vertex algebra
structure over $\C$ with ${\bf 1}$ as the vacuum vector and with
$$Y(f(z)u,x)=f(z+x)u(x),\ \ \ \
Y(f(z){\bf 1},x)=f(z+x)$$ for $u\in \g+\g^{1},\; f(z)\in \C((z))$.

Next, we show that $V_{\check{\g}_{p}}(\ell,0)$ is a vertex
$\C((z))$-algebra. Note that $V_{\check{\g}_{p}}(\ell,0)$ as a
vertex algebra over $\C$ is generated by
$\C((z))(\g+\g^{1})+\C((z)){\bf 1}$. For $f(z)\in \C((z)),\; u\in
\g+\g^{1}$, we have
$$Y(f(z)u,x)=f(z+x)Y(u,x),\ \ \ Y(f(z){\bf 1},x)=f(z+x).$$
Then it follows from \cite{li-tqva} that
$V_{\check{\g}_{p}}(\ell,0)$ is a vertex $\C((z))$-algebra.
\end{proof}

The following is a connection between restricted
$\hat{\g}_{p}$-modules of level $\ell$ and type zero
$V_{\check{\g}_{p}}(\ell,0)$-modules:

\bt{tmain} Let $\ell$ be a complex number. For any restricted
$\hat{\g}_{p}$-module $W$ of level $\ell$, there exists a unique
structure $Y_{W}$ of a type zero $V_{\check{\g}_{p}}(\ell,0)$-module
such that
$$Y_{W}(a,x)=a(x),\ \ Y_{W}(a^{1},x)=a^{1}(x)\ \
\mbox{ for }a\in \g.$$ On the other hand, let $(W,Y_{W})$ be a type
zero $V_{\check{\g}_{p}}(\ell,0)$-module. Then $W$ is a restricted
$\hat{\g}_{p}$-module of level $\ell$ with
$$a(x)=Y_{W}(a,x),\ \ \ a^{1}(x)=Y_{W}(a^{1},x)\ \ \ \mbox{ for
}a\in \g$$ and with $\k$ acting as scalar $\ell$.
 \et

\begin{proof} For the first assertion, the uniqueness is clear, since $\g+\g^{1}$
generates $V_{\check{\g}_{p}}(\ell,0)$ as a vertex
$\C((z))$-algebra. Set
$$U=\{a(x),\; a^{1}(x)\; |\; a\in \g\}.$$
{}From the defining relations of $\hat{\g}_{p}$, one sees that $U$
is a local subset of $\E(W)$. In view of Proposition \ref{pvta}, $U$
generates a vertex $\C((z))$-algebra $V_{W}=\C((x))\<U\>$ with $W$
as a faithful type zero module, where $f(z)\in \C((z))$ acts as
$f(x)$ on $V_{W}$. By Proposition \ref{pvacuum-module}, $V_{W}$ is a
vacuum $\check{\g}_{p}$-module of level $\ell$ with
\begin{eqnarray*}
a(x_{0})=Y_{\E}(a(x),x_{0}),\ \ \ \
a^{1}(x_{0})=Y_{\E}(a^{1}(x),x_{0})
\end{eqnarray*}
for $a\in \g$ and with $f(z)\in \C((z))$ acting as $f(x)$. {}From
the construction of $V_{\check{\g}_{p}}(\ell,0)$, there exists a
$\check{\g}_{p}$-module homomorphism $\psi$ from
$V_{\check{\g}_{p}}(\ell,0)$ to $V_{W}$ with $\psi({\bf 1})=1_{W}$.
Thus
\begin{eqnarray*}
&&\psi( u(n)v)=u(x)_{n}\psi(v)\ \ \ \mbox{ for }u\in \g+\g^{1},\;
n\in \Z,\; v\in V_{\check{\g}_{p}}(\ell,0),\\
&&\psi(f(z)v)=f(x)\psi(v)\ \ \ \mbox{ for }f(z)\in \C((z)).
\end{eqnarray*}
It follows that $\psi$ is a homomorphism of vertex
$\C((z))$-algebras. By Proposition \ref{pvta-module}, $W$ is a type
zero $V_{\check{\g}_{p}}(\ell,0)$-module.

On the other hand, let $(W,Y_{W})$ be a type zero
$V_{\check{\g}_{p}}(\ell,0)$-module. {}From definition, $(W,Y_{W})$
is a module for $V_{\check{\g}_{p}}(\ell,0)$ viewed as a vertex
algebra over $\C$. For $a,b\in \g$, we have
\begin{eqnarray*}
&&[Y(a^{1},x_{1}),Y(b^{1},x_{2})]\\
&=&Y(p(z)[a,b],x_{2})
x_{2}^{-1}\delta\left(\frac{x_{1}}{x_{2}}\right)
+\frac{1}{2}\<a,b\>\ell Y(p'(z){\bf 1},x_{2})
x_{2}^{-1}\delta\left(\frac{x_{1}}{x_{2}}\right)\\
&&\ \ +\<a,b\>\ell Y(p(z){\bf 1},x_{2})\frac{\partial}{\partial
x_{2}}x_{2}^{-1}\delta\left(\frac{x_{1}}{x_{2}}\right)\!.
\end{eqnarray*}
Thus
\begin{eqnarray*}
&&[Y_{W}(a^{1},x_{1}),Y_{W}(b^{1},x_{2})]\\
&=&Y_{W}(p(z)[a,b],x_{2})
x_{2}^{-1}\delta\left(\frac{x_{1}}{x_{2}}\right)
+\frac{1}{2}\<a,b\>\ell Y_{W}(p'(z){\bf 1},x_{2})
x_{2}^{-1}\delta\left(\frac{x_{1}}{x_{2}}\right)\\
&&\ \ +\<a,b\>\ell Y_{W}(p(z){\bf 1},x_{2})\frac{\partial}{\partial
x_{2}}x_{2}^{-1}\delta\left(\frac{x_{1}}{x_{2}}\right)\!\\
&=&p(x_{2})Y_{W}([a,b](x_{2})
x_{2}^{-1}\delta\left(\frac{x_{1}}{x_{2}}\right)
+\frac{1}{2}\<a,b\>\ell p'(x_{2})
x_{2}^{-1}\delta\left(\frac{x_{1}}{x_{2}}\right)\\
&&\ \ +\<a,b\>\ell p(x_{2})\frac{\partial}{\partial
x_{2}}x_{2}^{-1}\delta\left(\frac{x_{1}}{x_{2}}\right).
\end{eqnarray*}
The other relations are readily seen to hold. Therefore, $W$ is a
restricted $\hat{\g}_{p}$-module of level $\ell$.
\end{proof}

\end{document}